\documentclass[11pt]{article}
\usepackage[dvips]{graphics}
\usepackage{multirow}
\usepackage{array}
\usepackage{epsfig,rotating}
\usepackage{amssymb,amsmath}
\usepackage{latexsym}
\usepackage[usenames]{color}

\numberwithin{equation}{section}

\newcommand{\be}{\begin{equation}}
\newcommand{\ee}{\end{equation}}
\newcommand{\beaa}{\begin{eqnarray*}}
\newcommand{\eeaa}{\end{eqnarray*}}
\newcommand{\bea}{\begin{eqnarray}}
\newcommand{\eea}{\end{eqnarray}}
\newcommand{\lbl}{\label}

\newcommand{\bei}{\begin{itemize}}
\newcommand{\eei}{\end{itemize}}

\newcommand{\bd}{\bold}

\def\E{\mathrm{E}}
\def\Var{\mathrm{Var}}

\newtheorem{theorem}{ \noindent T{\footnotesize HEOREM}}
\newtheorem{prop}{ \noindent P{\footnotesize ROPOSITION}}
\newtheorem{lemma}{ \noindent L{\footnotesize EMMA}}
\newtheorem{coro}{ \noindent C{\footnotesize OROLLARY}}

\newcommand{\red}{\color{red}}

\begin{document}

\title{Empirical Distributions of Eigenvalues of Product Ensembles}
\author{Tiefeng Jiang$^{1}$
and  Yongcheng Qi$^2$\\
University of Minnesota}

\date{}
\maketitle

\footnotetext[1]{School of Statistics, University of Minnesota, 224 Church
Street, S. E., MN 55455, USA, jiang040@umn.edu. 
The research of Tiefeng Jiang was
supported in part by NSF Grant DMS-1209166 and DMS-1406279.}

\footnotetext[2]{Department of Mathematics and Statistics, University of Minnesota Duluth, MN 55812, USA, yqi@d.umn.edu.
The research of Yongcheng Qi was
supported in part by NSF Grant DMS-1005345.}

\begin{abstract}
\noindent  Assume a finite set of complex random variables form a determinantal point process, we obtain a theorem on the limit of the empirical distribution of these random variables. The result is applied to 
two types of  $n$ by $n$ random matrices as $n$ goes to infinity. The first one is the product of $m$ i.i.d. (complex) Ginibre ensembles, and the second one is the product of truncations of $m$ independent Haar unitary matrices with sizes $n_j\times n_j$ for $1\leq j \leq m$. Assuming $m$  depends on $n$, by using the special structures of the eigenvalues we developed, explicit limits of spectral distributions are obtained regardless of the speed of $m$ compared to $n$. For the product of $m$ Ginibre ensembles, as $m$ is fixed, the limiting distribution is known by various authors, e.g., G\"{o}tze and Tikhomirov (2010), Bordenave (2011), O'Rourke and Soshnikov (2011) and O'Rourke {\it et al}. (2014). Our results hold for  any $m\geq 1$ which may depend on $n$.
For the product of truncations of Haar-invariant unitary matrices, we show a rich feature of the limiting distribution as $n_j/n$'s vary.  In addition, some general results on arbitrary rotation-invariant determinantal point processes are also derived. Especially, we obtain an inequality for the fourth moment of linear statistics of  complex random variables forming a determinantal point process. This inequality is known for the complex Ginibre ensemble only [Hwang (1986)]. Our method is the determinantal point process rather than the contour integral by Hwang.

\end{abstract}

\noindent \textbf{Keywords:\/} non-symmetric random matrix, eigenvalue, empirical distribution,  determinantal point process.

\noindent\textbf{AMS 2000 Subject Classification:}
Primary: 15A52;\, Secondary: 60F99, 60G55, 62H10.\\


\newpage

\section{Introduction}
\setcounter{equation}{0}


In this paper we will study the limiting spectral laws of two types of random matrices. They are in the form of $\bd{X}_1\cdots \bd{X}_m$, which is called a product ensemble.  The first type is the product of $n\times n$ Ginibre ensembles, that is,  $\bd{X}_j,\, 1\leq j \leq m,$ are independent and identically distributed (i.i.d.) Ginibre ensembles; review that the $n\times n$ matrix $\bd{X}_1$ is referred to as a Ginibre ensemble if its $n^2$ entries are i.i.d. standard complex normal random variables. The second kind corresponds to that $\bd{X}_j,\, 1\leq j \leq m,$ are independent $n\times n$ matrices, each of which is a truncation of an Haar-invariant unitary matrix. We do not assume these matrices are of the same size.  To work on the two types of matrices, we derive a general result on complex random variables that form a determinantal point process. The limit of their empirical distribution can be obtained through the behavior of their radii only.


 After obtaining the general theorem (Theorem \ref{nonlinear}) mentioned above, we then investigate the structures of the eigenvalues of the two matrices (Lemmas \ref{Catroll} and \ref{new_loving}) by using a theory of the determinantal point processes. It is found that the absolute values of the eigenvalues are the product of i.i.d. Gamma-distributed random variables and the product of i.i.d. Beta-distributed random variables, respectively.

Using the theory, assuming $m$ depends on $n$, we obtain the limiting distributions of the eigenvalues of $\bd{X}_1\cdots \bd{X}_m$ for both cases as $n\to\infty$  by allowing $m$ to be fixed or go to infinity.
As $m$ does not depend on $n$ for the first case or $m=1$ for the second case, some knowledge about their limiting distributions are known. Here our results hold for any choice of  $m$. For the product of truncations of Haar unitary matrices with different sizes, the limiting distributions are very rich. 

The essential role in the derivation of our results is the determinantal point process $\{Z_1, \cdots, Z_n\}$. For the two product ensembles above, their kernels associated with the point process  are rotation-invariant. We then study it and obtain a general theory in Section \ref{Hu_page}. They may be useful in other occasions.

%

Before stating the main results, we need the following notation.

$\bullet$ Any function $g(z)$ of complex variable $z=x+iy$ should be interpreted as a bivariate function of $(x,y)$:  $g(z)=g(x, y)$.

$\bullet$ We write $\int_Ag(z)\,dz=\int_Ag(x,y)\,dxdy$ for any measurable set $A\subset \mathbb{C}.$

$\bullet$  $Unif(A)$ stands for the uniform distribution on a set $A$.

$\bullet$ For a sequence of random probability measures $\{\upsilon, \upsilon_n;\, n\geq 1\}$, we write
\bea\lbl{my_baby}
\upsilon_n \rightsquigarrow \upsilon\ \ \mbox{if \ P($\upsilon_n$  converges weakly to $\upsilon$ as $n\to\infty$)=1}.
\eea
When $\upsilon$ is a non-random probability measure generated by random variable $X$, we simply write $\upsilon_n \rightsquigarrow X$.
%
%
%
For complex variables $\{Z_j;\, 1\leq j \leq n\}$ mentioned above, we write
\bea\lbl{argument}
\Theta_j=\arg(Z_j)\in [0, 2\pi)\ \ \ \mbox{such that}\ \ \ Z_j=|Z_j|\cdot e^{i\Theta_j}
\eea
for each $j$. Let $Y_1, \cdots, Y_n$ be some given random variables, each of which may also rely on $n$. We omit the index $n$ for each $Y_j$ for clarity. Given a sequence of measurable functions $h_n(r), \, n\geq 1$, defined on $[0,\infty)$, set
\begin{equation}\label{mun}
\mu_n=\frac{1}{n}\sum^n_{j=1}\delta_{(\Theta_j, h_n(|Z_j|))}\ \ \mbox{and}\ \
\nu_n=\frac{1}{n}\sum^n_{j=1}\delta_{h_n(Y_j)}.
\end{equation}
 The empirical measure $\mu_n$ counts the frequency of the pairs of the angles and the radius of the $Z_j$'s. The measure $\nu_n$ counts the frequency of the $Y_j$'s. Roughly speaking, we can regard $Y_j$ as $|Z_j|$ for each $j$; see Lemma \ref{independence} latter.
 In (\ref{mun}), if $h_n$ is linear, that is
$h_n(r)=r/a_n$, where $\{a_n>0;\, n\geq 1\}$ is a sequence of numbers, we give special notation of the empirical measure of $Z_j$'s accordingly for this case by
\begin{equation}\label{mun*}
\mu_n^*=\frac{1}{n}\sum^n_{j=1}\delta_{Z_j/a_n}\ \ \mbox{and}\ \
\nu_n^*=\frac{1}{n}\sum^n_{j=1}\delta_{Y_j/a_n}.
\end{equation}


Review the notation ``$\rightsquigarrow$" in (\ref{my_baby}). The symbol $\mu_1\otimes\mu_2$ represents the product measure of two measures $\mu_1$ and $\mu_2$. Our general result is given as follows.

\begin{theorem}\label{nonlinear} Let $\varphi(x)\geq 0$ be a measurable  function defined on $[0, \infty).$ Assume the density  of $(Z_1, \cdots, Z_n)\in \mathbb{C}^n$ is proportional to  $\prod_{1\leq j < k
\leq n}|z_j-z_k|^2\cdot \prod_{j=1}^n\varphi(|z_j|)$. Let $Y_1, \cdots, Y_n$ be independent r.v.'s such that the density of $Y_j$ is proportional to  $y^{2j-1}\varphi(y)I(y\geq 0)$ for every $1\leq j\leq n.$ 
If $\{h_n\}$ are  measurable functions such that $\nu_n\rightsquigarrow\nu$ for some probability measure $\nu$, then $\mu_n \rightsquigarrow \mu$ with
$\mu=Unif[0, 2\pi]\otimes\nu$ .
Taking $h_n(r)=r/a_n$, the conclusion still holds if ``$(\mu_n, \nu_n, \mu, \nu)$" is replaced by  ``$(\mu_n^*, \nu_n^*, \mu^*, \nu^*)$" where $\mu^*$ is the distribution of $Re^{i\Theta}$ with $(\Theta, R)$ having the law of  $Unif[0, 2\pi]\otimes\nu^*$.

\end{theorem}

Next we apply this theorem to two types of product ensembles. The first one is the product of $m$ Ginibre ensembles given in Section \ref{1212}. The second one is the product of $m$ truncated unitary matrices presented in Section \ref{2323}. For the first one, no results are known  as $m$ depends on $n$ and $m\to\infty$. We actually will give a universal result regardless of the speed of $m$ relative to $n.$ For the second product ensemble, it is not clear whether the empirical distributions of eigenvalues converge or not in the previous literature. We obtain the limiting  laws, which have an interesting feature: the limiting law are  very different when the sizes of truncations vary.


\subsection{Product of Ginibre Ensembles}\lbl{1212}

Given an integer $m\geq 1$. Assume $\bd{X}_1, \cdots, \bd{X}_m$ are i.i.d. $n\times n$ random matrices and the $n^2$ entries of $\bd{X}_1$ are i.i.d. with the standard complex normal distribution $\mathbb{C}N(0, 1).$ Let $Z_1, \cdots, Z_n$ be the eigenvalues of the product $\prod_{j=1}^m\bd{X}_j$. It is known that their joint density function is
\bea\label{up_stairs}
C\prod_{1\leq j <l \leq n}|z_j - z_l|^2\prod_{j=1}^nw_{m}(|z_j|)
\eea
where $C$ is a normalizing constant and $w_m(z)$  has a recursive formula given by
$w_1(z)=\exp(-|z|^2)$ and
\beaa
w_{m}(z)=2\pi\int^\infty_0w_{m-1}\Big(\frac{z}{r}\Big)\exp(-r^2)\frac{dr}{r}
\eeaa
for all integer $m\geq 2$; see, e.g., Akemann and Burda (2012). The function $w_m(z)$ also has a representation in terms of the so-called Meijer G-function; see the previous reference.

Through investigating the limit of the kernel of a determinantal point process, Burda {\it et al}. (2010) and  Burda (2013) showed that  the empirical distribution of $Z_j/n^{m/2},\, 1\leq j \leq n$, in the sense of mean value, converges to a distribution with density $\frac{1}{m\pi}|z|^{\frac{2}{m}-2}$ for $|z|\leq 1.$ Under the condition that the entries of $\bd{X}_1$ are i.i.d. random variables with a certain moment condition, G\"{o}tze and Tikhomirov (2010) prove  the above result in the sense of mean value. Bordenave (2011), O'Rourke and Soshnikov (2011) and O'Rourke {\it et al}. (2014) further generalize this result to the almost sure convergence. Our result next gives a weak convergence of the eigenvalues $Z_j$'s by allowing  $m$ to depend on $n$, and the result holds regardless of the speed of $m$ relative to $n$. Review (\ref{argument}).


\begin{theorem}\label{Theorem_product} Let $\{m_n\geq 1;\, n\geq 1\}$ be an arbitrary sequence of integers.
Define
\[
\mu_n=\frac1n\sum^n_{j=1}\delta_{(\Theta_j, \frac{1}{n}\,|Z_j|^{2/m_n})}.
\]
Then $\mu_n \rightsquigarrow Unif\big([0, 2\pi)\otimes [0,1]\big)$ as $n\to\infty$.
\end{theorem}
Theorem \ref{Theorem_product} implies that the angle and the length of a randomly picked pair $(\Theta_j, \frac{1}{n}\,|Z_j|^{2/m_n})$ are asymptotically independent. Take  $m_n= m$ for all $n\geq 1$. By the continuous mapping theorem,  the above conclusion implies that, with probability one, the empirical distribution of $\frac{1}{n^{m/2}}Z_j=\frac{1}{n^{m/2}}\,|Z_j|\,e^{i\Theta_j}$, $1\leq j \leq n$, converges weakly to the distribution of $R^{m/2}e^{i\Theta}$, where $(\Theta, R)$ follows the law $Unif\big([0, 2\pi)\times [0,1]\big)$. Easily, $Z:=R^{m/2}e^{i\Theta}\in \mathbb{C}$ has density $\frac{1}{m\pi}|z|^{\frac{2}{m}-2}$ for $|z|\leq 1.$ This yields the conclusion mentioned  before Theorem \ref{Theorem_product}.  In particular, taking  $m_n= 1$ for all $n\geq 1$, we have
\bea\lbl{circular_law}
\frac1n\sum^n_{j=1}\delta_{Z_j/\sqrt{n}} \rightsquigarrow Unif\{|z|\leq 1\}.
\eea
This gives the classical circular law.  For the universality of (\ref{circular_law}), where the entries of $\bd{X}_1$ are not necessarily Gaussian, one can check, for instance, Girko (1984), Bai (1997), Tao and Vu (2010) or Bordenave and Chafa\"{\i} (2012).

 %

The proof of Theorem \ref{Theorem_product} is based on Theorem \ref{nonlinear} and a recent result by  Jiang and Qi (2017) in which the exact distributions of $|Z_j|$ are shown to be the products of independent Gamma-distributed random variables (see Lemma \ref{Catroll}).

\subsection{Products of Truncated Unitary Matrices}\lbl{2323}

Let $m$, $n$ and $\{l_j; 1\leq j \leq m\}$ be positive integers. Set $n_j=l_j+n$ for $1\le j\le m.$ Suppose $\{U_j;\,1\le j\le m\}$ are independent Haar-invariant unitary matrices where $U_j$ is $n_j\times n_j$ for each $j$. Let $\mathbf{X}_j$ be the upper-left $n\times n$ sub-matrix of $U_j$. Consider  the product matrix $\mathbf{X}^{(m)}=\mathbf{X}_m\mathbf{X}_{m-1}\cdots\mathbf{X}_1$. We write this way instead of the product of the matrices in the reverse order is simply for brevity of notation below.

The joint density function for the eigenvalues $Z_1, \cdots, Z_n$ of  $\mathbf{X}^{(m)}$ is derived by
Akemann {\it et al}. (2014):
\begin{equation}\label{truncation-density}
p(z_1, \cdots, z_n)=C\prod_{1\le j<k\le n}|z_j-z_k|^2\prod^n_{j=1}w_m^{(l_1, \cdots, l_m)}(|z_j|)
\end{equation}
for all $z_j$'s with $\max_{1\leq j \leq m}|z_j|<1$, where $C=\frac{1}{n!}\prod^m_{j=0}\prod^{n-1}_{l=0}{{l_j+l}\choose{l}}^{-1}$, and $w_m^{(l_1, \cdots, l_m)}$ can be recursively obtained by
\bea\lbl{salt_sugar}
w_m^{(l_1, \cdots, l_m)}(s)=2\pi\int^1_0w_1^{(l_m)}(u)w_{m-1}^{(l_1, \cdots, l_{m-1})}(\frac{s}{u})\frac{du}{u}, ~~~s\in [0,1),
\eea
with initial $w_1^{(l)}(s)=(l/\pi)(1-s^2)^{l-1}I(0\le s<1)$. The function  $w_m^{(l_1, \cdots, l_m)}(s)$ can be expressed in terms of Meijer G-functions. One can see Appendix C from  Akemann {\it et al}. (2014) for details. The density in (\ref{truncation-density}) for the case $m=1$ is obtained by \.{Z}yczkowski and Sommers (2000).

Although the density of the eigenvalues is given in (\ref{truncation-density}), no limit law of the empirical distribution of  $Z_1, \cdots, Z_n$ is known, even heuristic results. We now consider the problem.  Assume  $m$ depends on $n$  and
$1<n<\min_{1\le j\le m}n_j$.  For convenience, we assume
$n_1, \cdots, n_m$ are functions of $n$ and all limits are taken as $n\to\infty$ unless otherwise specified. The limiting spectral distribution actually depends on the limit of functions $F_n(x)$'s defined below.

Let $\{\gamma_n;\, ~n\ge 1\}$ be a sequence of positive numbers.  Define
\begin{equation*}
g_n(x)=\prod^{m}_{j=1}\frac{nx}{nx+l_j}, ~~~~x\in [0, 1]
\end{equation*}
and
\begin{equation}\label{fnx}
F_n(x)=\Big(\frac{g_n(x)}{g_n(1)}\Big)^{1/\gamma_n}=\Big(\prod^{m}_{j=1}\frac{n_jx}{nx+l_j}\Big)^{1/\gamma_n}, ~~~~~x\in [0, 1].
\end{equation}
Note that $F_n(x)$ is continuous and strictly increasing on $[0,1]$, $F_n(0)=0$ and $F_n(1)=1$.
We will assume that $F_n(x)$ has a limit $F(x)$  defined on $(0, 1)$ such that
\bea
& &  \mbox{$F(x)$ is continuous and strictly increasing over (0,1)},\nonumber\\
& & \mbox{$\lim_{x\downarrow 0}F(x)=0$\ \ \mbox{and}\ \ $\lim_{x\uparrow 1}F(x)=1$}
\ \ \ \ \ \ \ \  \lbl{C_1_condition}
\eea
or a limit $F(x)$ defined on $(0, 1]$ satisfying
\bea
  \mbox{$F(x)=1$ for all $x\in (0,1]$}.\lbl{C_2_condition}
\eea

Recall the notation ``$\rightsquigarrow$" in (\ref{my_baby}) and ``$Unif(A)$" standing for the uniform distribution on a set $A$. Write $z=re^{i\theta}$. Evidently, $z\in \mathbb{C}$ and $(\theta, r)\in [0, 2\pi) \times [0, \infty)$ are one-to-one correspondent. Although the empirical distributions considered in the following  are targeted as functions of  $(\theta, r)$, we will characterize their limits  in terms of complex distributions since the arc law and the circular law etc are easily understood.

 Review (\ref{C_1_condition}). In the following result, we denote $f^*(x)=\frac{d}{dx}F^{-1}(x)$ for $0<x<1$ if the derivative exists.
\begin{theorem}\label{pie}
Assume there exists a sequence of numbers $\{\gamma_n\}$ with $\gamma_n\ge 1$  such that
$\lim_{n\to\infty}F_n(x)=F(x), ~x\in (0,1)$,
for some function $F(x)$ defined on $(0,1)$. Let $Z_1, \cdots, Z_n$ be the eigenvalues of $\bd{X}^{(m)}$, $\mu_n$ be as in \eqref{mun} with $h_n(r)=(r^2/b_n)^{1/\gamma_n}$ and  $b_n=\prod^{m}_{j=1}\frac{n}{n_j}$.


\noindent{\textbf{(a)}.} If (\ref{C_1_condition}) holds and  $f^*(x)$
 exists for $x\in (0,1),$ then
$\mu_n \rightsquigarrow (\Theta, R)$ such that $Z:=Re^{i\Theta}$ has density $\frac{1}{2\pi|z|}f^*(|z|)$ for $0<|z|< 1$.



\noindent{\textbf{(b)}.} If (\ref{C_2_condition}) holds,  then
$\mu_n \rightsquigarrow (\Theta, R)$ such that $Z=Re^{i\Theta}$ has the law $Unif\{|z|=1\}$.

\noindent{\textbf{(c)}.}  If $\gamma_n=2$, (\ref{C_1_condition}) holds and  $f^*(x)$ exists for $x\in (0,1)$, then
$\mu_n \rightsquigarrow (\Theta, R)$ such that $Z=Re^{i\Theta}$ has density $\frac{1}{2\pi|z|}f^*(|z|)$ for $0<|z|< 1$. Further, let $\mu_n^*$ be as in \eqref{mun*} with $a_n=(\prod^{m}_{j=1}\frac{n}{n_j})^{1/2}$. If $\gamma_n=2$  and (\ref{C_2_condition}) holds, then
$\mu_n^* \rightsquigarrow \mu^*$
with $\mu^*= Unif\{|z|=1\}$.
\end{theorem}

From the proof of Theorem \ref{pie}, it is actually seen that the condition ``$F^{-1}(x)$ is differentiable" is not necessary in (a). The general conclusion is that   $\mu_n \rightsquigarrow\mu$ where $\mu$ is the product measure of $Unif[0, 2\pi]$ and the probability measure on $[0, \infty)$ with cumulative distribution function $F^{-1}(x)$. We write the way in (a) to avoid a lengthy statement. In particular, (a) is general enough to our applications next.
%

The values of $n_j$'s in Theorem \ref{pie} can be very different. Now let us play them and find out their limiting distributions. The first one below is on the size $n_j$'s that are at the same scale and $m$ is fixed.
\begin{coro}\label{cor1}  Assume that $m\geq 1$ is an integer and $m$ does not depend on $n$ and that $\lim_{n\to\infty}\frac{n}{n_j}=\alpha_j\in [0,1]$ for $1\le j\le m$.  Assume $Z_1, \cdots, Z_n$ are the eigenvalues of $\bd{X}^{(m)}$. Let $a_n=(\prod^m_{j=1}\frac{n}{n_j})^{1/2}$ and $\mu_n^*$ be as in (\ref{mun*}).

\noindent(1). If $\alpha_1=\cdots=\alpha_m=1$, then $\mu^*_n\rightsquigarrow Unif\{|z|=1\}$.

\noindent(2). If $\alpha_1=\cdots=\alpha_m=\alpha\in [0,1)$, then $\mu^*_n\rightsquigarrow\mu^*$ with density $\frac{1}{2\pi|z|}f^*(|z|)$ on $\{0<|z|< 1\}$, where $f^*(x)=2(1-\alpha)m^{-1}x^{(2/m)-1}(1-\alpha x^{2/m})^{-2}$ for $x\in [0, 1].$

\noindent(3). If $\alpha_j < 1$ for some $1\leq j \leq m$, then $\mu^*_n \rightsquigarrow\mu^*$ with  density $\frac{1}{2\pi|z|}f^*(|z|)$ on the set $\{0<|z|< 1\}$ where  $f^*(x)=\frac{d}{dx}F^{-1}(x)$.
\end{coro}

Trivially, part (2) in the above corollary is a special case of part (3). We single it out since $f^*(x)$ has an explicit expression. Picking $m=1$ in (1) of Corollary \ref{cor1}, we know  that, with probability one,  $\mu^*_n\rightsquigarrow Unif\{|z|=1\}$, where $\mu_n^*=\frac{1}{n}\sum^n_{j=1}\delta_{Z_j/a_n}$
with $a_n=(n/n_1)^{1/2}\to 1$ as $n\to\infty.$ It implies that $\frac{1}{n}\sum^n_{j=1}\delta_{Z_j}$ converges weakly to $Unif\{|z|=1\}$. This conclusion is obtained by Dong {\it et al}. (2012). Taking $m=1$ and $\alpha_1\in (0,1)$ in (2) of Corollary \ref{cor1}, we get a result by Petz and R\'{e}ffy (2005):
$\frac{1}{n}\sum_{j=1}^n\delta_{Z_j} $ converges weakly to a probability measure with
density $f(z)=\frac{\alpha_1^{-1}-1}{\pi (1-|z|^2)^2}$ for $|z|\leq \sqrt{\alpha_1}$ (noticing the scaling in $\mu^*_n$ is different from $\frac{1}{n}\sum_{j=1}^n\delta_{Z_j}$).

Pick $\alpha=0$ from (2) of Corollary \ref{cor1}, the limiting density becomes $\frac{1}{m\pi}|z|^{\frac{2}{m}-2}$ for $|z|\leq 1$,
which is exactly the same as that of the product of Ginibre ensembles; see the paragraph above Theorem \ref{Theorem_product}. This is not a coincidence. In fact, Jiang (2009) show that the $n\times n$ submatrix  $\bd{X}_1$ of the $n_1\times n_1$ matrix $U_1$ can be approximated by a Ginibre ensemble as $n=o(\sqrt{n_1})$ in the variation norm. Similar conclusion also holds for Haar-invariant orthogonal matrices (Jiang, 2006).

If  $m$ depends on $n$ and $m\to \infty$, and $\{n_j;\, 1\leq j \leq m\}$ are almost sitting on a curve, what is the corresponding limit appearing in Theorem \ref{pie}? To answer the question, assume there exists a continuous function
$q(x)$ defined over $[0,1]$ satisfying  $0<q(x)<1$ for $0<x<1$ and
\begin{equation}\label{case22}
\lim_{n\to\infty}\frac{1}{m}\sum^{m}_{j=1}\Big|\frac{n}{n_j}-q\Big(\frac{j}{m}\Big)\Big|=0.
\end{equation}
Define
\bea
& & F(x)=x\exp\Big(-\int^1_0\log\big(1-q(t)(1-x)\big)dt\Big); \lbl{twin_star}\\
& & f(x)=\frac{F(x)}{x}\int^1_0\frac{1-q(t)}{1-q(t)(1-x)}dt \lbl{Sly}
\eea
for $0< x\le 1$ and $f(0)=F(0)=0$.

\begin{coro}\label{cor2} Assume  $m$ depends on $n$, $m\to \infty$ and (\ref{case22}) holds. Let $F(x)$ and $f(x)$ be as in (\ref{twin_star}) and (\ref{Sly}), respectively. Let $Z_1, \cdots, Z_n$ be the eigenvalues of $\bd{X}^{(m)}$. Set $b_n=\prod^{m}_{j=1}(n/n_j)$ and
\[
\mu_n=\frac{1}{n}\sum^n_{j=1}\delta_{(\Theta_j, (|Z_j|^2/b_n)^{1/m})}.
\]
Then
$\mu_n \rightsquigarrow (\Theta, R)$ and $Z=Re^{i\Theta}$ has density $\frac{1}{2\pi|z|\,f(F^{-1}(|z|))}$ for $0<|z|<1.$
%
\end{coro}


For an Haar-invariant unitary matrix, the empirical distribution of its eigenvalues is asymptotically the arc law $Unif\{|z|=1\}$, see, for example, Diaconis and  Shahshahani (1994) and  Diaconis and Evans (2001). If $n$ is very close to $n_j$ for each $j$ in Theorem \ref{pie}, that is, the truncated sub-matrix with size $n\times n$ of  $\bd{U}_j$ is almost the same as $\bd{U}_j$ for each $j$, do we always expect the arc law $Unif\{|z|=1\}$? The answer is no and, as a matter of fact, it depends on the sum of $l_j=n_j-n$ for $1\leq j \leq m.$

\begin{coro}\label{cor3} Let  $m$ depend on $n$, $m\to \infty$
and $\lim_{n\to\infty}\max_{1\le j\le m}|\frac{n}{n_j}-1|=0$.
Let $\lim_{n\to\infty}\frac1n\sum^{m}_{j=1}l_j=\beta\in [0,\infty]$. Assume $Z_1, \cdots, Z_n$ are the eigenvalues of $\bd{X}^{(m)}$. Let $\mu_n^*$ be as in (\ref{mun*}) with $a_n=(\prod^{m}_{j=1}\frac{n}{n_j})^{1/2}$.

\noindent (a).  If $\beta=0$,
then $\mu^*_n \rightsquigarrow Unif\{|z|=1\}$.

Let  $\mu_n$ be as in \eqref{mun} with $h_n(r)=(r^2/b_n)^{1/\gamma_n}$ and    $b_n=\prod^{m}_{j=1}\frac{n}{n_j}$.

\noindent (b).  If $\beta\in (0,\infty)$, with $\gamma_n=2$ we have
$\mu_n\rightsquigarrow\mu$  with density $\frac{\beta}{\pi|z|^2(\beta-2\log |z|)^2}$ for $0<|z|<1$.
\newline
\noindent (c). If $\beta=\infty$, with $\gamma_n=\frac1n\sum^{m}_{j=1}l_j$ we have $\mu_n \rightsquigarrow\mu$ where $\mu$ has density $\frac{1}{2\pi |z|^{2}(1-\log |z|)^2}$ for $0<|z|<1$.
\end{coro}


Finally, we work on the case that $n$ is much smaller than $n_j$'s.

\begin{coro}\label{cor4}  Let $m$ depend on $n$, $m\to \infty$ and $\max_{1\le j\le m}\frac{n}{n_j}=0$ as $n\to\infty$. Review $Z_1, \cdots, Z_n$ are the eigenvalues of $\bd{X}^{(m)}$. Set $b_n=\prod^{m}_{j=1}(n/n_j)$ and
\[
\mu_n=\frac{1}{n}\sum^n_{j=1}\delta_{(\Theta_j, (|Z_j|^2/b_n)^{1/m})}.
\]
Then
$\mu_n  \rightsquigarrow (\Theta, R)$ and $Z=\sqrt{R}e^{i\Theta}$ follows $Unif\{ |z|\leq 1\}$, that is, the circular law.
\end{coro}

Picking $m=1$, since $Z_j=|Z_j|e^{i\Theta_j}$, by the continuous mapping theorem, we get that, with probability one,  $\frac{1}{n}\sum_{j=1}^n\delta_{(n_j/n)^{1/2}Z_j}$ converges weakly to the circular law $Unif\{ |z|\leq 1\}$. This result is found and proved by Dong {\it et al}. (2012).

\subsection{Structures of Determinantal Point Processes on Complex Plane}\lbl{Hu_page}

In this section we state our results on  rotation-invariant determinantal point processes on complex plane; see the set-up in Lemma \ref{independence} below. The normalizing constant of their joint density function, moments and the structures of the two product matrices aforementioned are obtained.

Let $\{Z_1, \cdots, Z_n\}$ be $n$ complex-valued random variables. Let $K(z, w): \mathbb{C}^2\to \mathbb{C}$ with $\overline{K(z,w)}=K(w,z)$ for all $(z, w)\in \mathbb{C}^2.$ Let $\nu$ be a Borel  measure on $\mathbb{C}.$ We say $\{Z_1, \cdots, Z_n\}$ forms a determinantal point process with kernel $K(z,w)$ and background measure $\nu$ if the density function of $\{Z_1, \cdots, Z_k\}$ is given by
\bea\lbl{Space_length}
f_n(z_1, \cdots, z_k):=\frac{(n-k)!}{n!}\,\mbox{det}\big(K(z_j, z_l)\big)_{1\leq j, l\leq k},\ \ (z_1, \cdots, z_k) \in \mathbb{C}^k,
\eea
with respect to the product measure $\nu^{\otimes k}$ for all $1\leq k \leq n.$ The determinantal point process given above is a special case of a general definition in which the space $\mathbb{C}$ can be a arbitrary metric space. The definition here is good enough for our discussion. For general case, one can see, for example, Soshnikov (2000), Johansson (2005) or  Hough {\it et al}. (2009) for a reference.

Let $\varphi(x)\geq 0$ be a measurable  function defined on $[0, \infty)$ with   $0<\int_0^{\infty}y^{2j-1}\varphi(y)\,dy<\infty$ for each $1\leq j \leq n.$ Define
\begin{equation}\label{densityofYj}
p_j(y)=\frac{y^{2j-1}\varphi(y)I(y\geq 0)}{\int^\infty_0y^{2j-1}\varphi(y)dy}
\end{equation}
for $1\le j\le n$.  Define
\bea\lbl{jilin}
P_n(y)=\frac{1}{n}\sum^n_{j=1}p_j(y), \ \ y\geq 0.
\eea

Now we start a series of results on the determinantal point processes. The following is a special case of Theorem 1.2 from  Chafa\"{\i} and P\'{e}ch\'{e} (2014). It is another version of  Theorem 4.7.1 from Hough {\it et al}. (2009).
\begin{lemma}\lbl{independence} Let $\varphi(x)\geq 0$ be a measurable  function defined on $[0, \infty).$ Let $f(z_1, \cdots, z_n)$ be the probability density function of $(Z_1, \cdots, Z_n)\in \mathbb{C}^n$ such that it is  proportional to  $\prod_{1\leq j < k
\leq n}|z_j-z_k|^2\cdot \prod_{j=1}^n\varphi(|z_j|)$. Let $Y_1, \cdots, Y_n$ be independent r.v.'s such that the density of $Y_j$ is proportional to  $y^{2j-1}\varphi(y)I(y\geq 0)$ for each $1\leq j \leq n.$ Then,  $g(|Z_1|, \cdots, |Z_n|)$ and $g(Y_1, \cdots, Y_n)$ have the same distribution for any symmetric function $g(y_1, \cdots, y_n)$.
\end{lemma}

The next result is mostly known. Our contribution is that we are able to evaluate $C$ for any $\varphi(x)$, where $f(z_1, \cdots, z_n)=C\cdot\prod_{1\leq j < k
\leq n}|z_j-z_k|^2\cdot \prod_{j=1}^n\varphi(|z_j|)$ is as in Lemma \ref{independence}.

\begin{prop}\lbl{Laplace} Let $\varphi$ and $f$ be as in Lemma \ref{independence}.
Set $c_k=2\pi\int_0^{\infty}x^{2k+1}\varphi(x)\,dx$ for all $k=0,\cdots, n-1.$ Then, $C^{-1}=n!c_0c_1\cdots c_{n-1}$ and $(Z_1, \cdots, Z_n)$ forms a determinantal point process with background measure $\varphi(|z|)\,dz$ and kernel $K(z, w)=\sum_{k=0}^{n-1}\frac{1}{c_k}(z\bar{w})^k.$
\end{prop}

The next result gives an estimate of the fourth moment of the sum of a function of $Z_j$'s, where the $Z_j$'s forms a determinantal point process. Hwang (1986) obtains a similar result for the special case of the complex  Ginibre ensemble with $\varphi(z)=e^{-|z|^2}$. In particular, we do not assume any differentiability  of $\varphi(x).$

\begin{prop}\lbl{Aha_you} Let $Z_1, \cdots, Z_n$ and $\varphi(x)$ be  as in Proposition \ref{Laplace} with $\mu(dz)=\varphi(|z|)dz$.
Then,  for  any measurable function $h(z):\, \mathbb{C}\to \mathbb{R}$ with
$\sup_{z\in \mathbb{C}}|h(z)|\le 1$, we have
\beaa
\E\Big[\sum^n_{j=1}(h(Z_j)-\E h(Z_j))\Big]^4\le K n^2
\eeaa
for all $n\ge 1$, where $K$ is a constant not depending on $n$, $\varphi(z)$ or $h(z).$
\end{prop}
The essential of the proof of Proposition \ref{Aha_you} is the estimate of $\E\prod_{j=1}^4(h(Z_j)-\E h(Z_j))$. It is carried out by using (\ref{Space_length}) repeatedly. Our proof is different from the analysis of a contour integral by Hwang (1986), which seems to fit the Gaussian kernel $\varphi(z)=e^{-|z|^2}$ only.

 Let $U_1, \cdots, U_n$ be independent and real-valued random variables with $\mbox{Var}(U_j)=1$ for each $j$ and $C:=\sup_{1\leq j \leq n}E(U_j^4)<\infty.$ It is easy to check that $\E[\sum^n_{j=1}(U_j-\E U_j)]^4\le K n^2$, where $K$ is a constant depending on $C$ but not depending on $n.$ This suggests that, although the $Z_j$'s from Proposition \ref{Aha_you} are correlated each other, they are weakly correlated.

We will use Proposition \ref{Aha_you}, the Markov inequality and the Borel-Cantelli lemma to prove the almost sure convergence stated in Theorem~\ref{nonlinear}.

\begin{prop}\label{lem1} Let $Z_i$'s and $\varphi(x)$ be  as in Lemma \ref{independence}. Then the following hold.\\
(i) Let $P_n(\cdot)$ be as in (\ref{jilin}). Then, for any bounded measurable function $h(z)$,
\[
\E h(Z_1)=\frac{1}{2\pi}\int^\infty_{0}\Big(\int^{2\pi}_0h(re^{i\theta})d\theta \Big)P_n(r)dr.
\]
(ii) Let $\Theta_1$ be as in (\ref{argument}). Then, $|Z_1|$ has density $P_n(r)$, $\Theta_1 \sim Unif[0, 2\pi]$ and the two are independent.  Consequently, for any bounded measurable function $g(r, \theta)$,
\bea\lbl{Coke_pepsi}
\E g(\Theta_1, |Z_1|)= \frac{1}{2\pi}\int^\infty_{0}\Big(\int^{2\pi}_0g(\theta, r)d\theta \Big)P_n(r)dr.
\eea
\end{prop}
The following result reveals the structure of the eigenvalues of the product of Ginibre ensembles. The key is the Gamma distribution. We will switch the use of the notation of ``$r$"  in the next two lemmas, which will serve as an index instead of the radius of a complex number used earlier.

\begin{lemma}\lbl{Catroll} Let $Z_1,\cdots, Z_n$ have density as in  (\ref{up_stairs}). Let $\{s_{j, r},\,1\leq j\leq  n,\, 1\le r\le m\}$ be independent r.v.'s and $s_{j, r}$ have the Gamma  density $y^{j-1}e^{-y}I(y>0)/(j-1)!$ for each $j$ and $r$.
Then $g(|Z_1|^2, \cdots, |Z_n|^2)$ and $g(\prod_{r=1}^ms_{1,r}, \cdots, \prod_{r=1}^ms_{n,r})$ have the same distribution for any symmetric function $g(t_1, \cdots, t_n)$.
\end{lemma}

Recall the beta function
\bea\label{beta-gamma}
B(a, b)=\int^1_0s^{a-1}(1-s)^{b-1}\,ds=\frac{\Gamma(a)\Gamma(b)}{\Gamma(a+b)},\ \ a>0,\ b>0.
\eea

The following lemma describes the structure of the eigenvalues of the product of truncations of Haar-invariant unitary matrices. It has the same setting as Lemma \ref{Catroll} with ``Gamma distribution" replaced by ``Beta distribution".
\begin{lemma}\lbl{new_loving} Let $Z_1,\cdots, Z_n$ have density as in (\ref{truncation-density}). Let $\{s_{j, r},\,1\leq j\leq n,\, 1\le r\le m\}$ be independent r.v.'s and $s_{j, r}$ have the Beta  density $\frac{1}{B(j,l_r)}y^{j-1}(1-y)^{l_r-1}I(0\leq y\leq 1)$ for each $j$ and $r$.
Then $g(|Z_1|^2, \cdots, |Z_n|^2)$ and $g(\prod_{r=1}^ms_{1,r}, \cdots, \prod_{r=1}^ms_{n,r})$ have the same distribution for any symmetric function $g(t_1, \cdots, t_n)$.
\end{lemma}

\medskip

\noindent\textbf{Comments}. In this paper we study the empirical distribution of the eigenvalues of the product of $m$ Ginibre ensembles as well as the product of truncations of $m$ Haar-invariant unitary matrices.
 Now we make some remarks and state certain problems for  future.

1. There are other type of studies on the product of random matrices in literature. The size of each matrix is assumed to be  fixed and conclusions are obtained by letting the number of matrices go to infinity. Two typical interests of the product matrices are their norms and entries; see, for example, Furstenberg and Kesten (1960) and Mukherjea (2000).

2. In this paper we study two kinds of product matrices: the product of Ginibre ensembles and that of truncated Haar unitary matrices. Notice the Ginibre ensemble and  truncated Haar unitary matrices are of the Haar-invariant property. We believe that the same method (Theorem \ref{nonlinear}) can be used to derive the spectral limits of the products of other independent Haar-invariant matrices. The key is the explicit formula of $\varphi(x)$ and the verification of ``$\nu_n\rightsquigarrow\nu$" as stated in Theorem \ref{nonlinear}.

3. The universality of Theorem \ref{Theorem_product} is an interesting problem. Namely, replacing the normal entries in the Ginibre ensemble with i.i.d. non-Gaussian random variables, does Theorem \ref{Theorem_product} still hold?  In fact, Bordenave (2011), O'Rourke and Soshnikov (2011) and O'Rourke {\it et al}. (2014) show that this is true for fixed $m$. We expect  Theorem \ref{Theorem_product} to hold for non-Gaussian entries and arbitrary $m$ which may depend on $n$.

4. There are three Haar-invariant matrices generating the Haar measure of the classical compact groups: Haar-invariant orthogonal, unitary and symplectic matrices; see, for example, Jiang (2009, 2010). Similar to Theorem \ref{pie} one can work on the same limiting problems for the orthogonal and symplectic matrices.

5. If we change the square matrices $\bd{X}_j$ in Theorem \ref{Theorem_product}  to rectangular matrices and keep the Gaussian entries of each matrix, that is, $\bd{X}_j$ is $n_j\times n_{j+1}$ with $n_{m+1}=n_1$, it will be interesting to see the corresponding result. The limiting distribution will have a rich feature as the ratio $n_j/n_{j+1}$ fluctuates for each $j$.

6. Let $Z_1, \cdots, Z_n$ be the eigenvalues of the product $\prod_{j=1}^m\bd{X}_j$, where $\bd{X}_1, \cdots, \bd{X}_m$ are i.i.d. Ginibre ensembles. The spectral distribution of $Z_1, \cdots, Z_n$ is well understood through Theorem \ref{Theorem_product}. The transitional phenomenon of the spectral radius $\max_{1\leq i \leq m}|Z_j|$ is obtained by Jiang and Qi (2017). It is classified by $c:=\lim \frac{m}{n}$ with $c=0$, $c\in (0, \infty)$ and $c=\infty.$ The spectral radius of $\bd{X}^{(m)}$ in Theorem \ref{pie} with $m=1$ is also investigated in the same paper. With the help of Lemma \ref{new_loving}, the spectral radius of $\bd{X}^{(m)}$ for arbitrary $m$ can be done similarly.

The rest of the paper is organized as follows. In Section \ref{new_do}, we prove the results stated in Section \ref{Hu_page} that serve as technical tools to prove the main results. We then prove  Theorems \ref{nonlinear}-\ref{pie} and all of the corollaries in Sections \ref{yes_right}-\ref{god_bite}.

\section{Proofs}

In this section we will prove the main results stated in the Introduction. We first prove those
 in Section \ref{Hu_page} since they serve as
tools to derive the main limit theorems. Theorems \ref{nonlinear}-\ref{pie} and all corollaries are proved one by one afterwards.

\subsection{Proofs of Propositions \ref{Laplace}-\ref{lem1}, Lemmas \ref{Catroll} and \ref{new_loving}}\lbl{new_do}

\noindent\textbf{Proof of Proposition \ref{Laplace}.}
For $1\leq j \leq n$, let $Z_j=R_je^{i\Theta_j}$ with $R_j\geq 0$ and $\Theta_j\in [0, 2\pi).$ The Jacobian is obviously equal to $r_1\cdots r_n.$ Thus the joint density function of $(R_1, \cdots, R_n, \Theta_1, \cdots, \Theta_n)$ is given by
\bea\lbl{mean}
C\prod_{1\leq j < k
\leq n}|r_je^{i\theta_j}-r_ke^{i\theta_k}|^2\cdot \prod_{j=1}^n\big(r_j\varphi(r_j)\big).
\eea
Set
\beaa
M_n=
\begin{pmatrix}
1 & 1 & \cdots & 1\\
r_1e^{i\theta_1} & r_2e^{i\theta_2} & \cdots & r_ne^{i\theta_n}\\
\vdots\\
r_1^{n-1}e^{i(n-1)\theta_1} & r_2^{n-1}e^{i(n-1)\theta_2} & \cdots & r_n^{n-1}e^{i(n-1)\theta_n}
\end{pmatrix}
.
\eeaa
By the formula of the Vandermonde determinant, the first product in (\ref{mean}) is equal to $(\mbox{det}(M_n))^2.$ It follows that
\beaa
\prod_{1\leq j < k
\leq n}|r_je^{i\theta_j}-r_ke^{i\theta_k}|^2 = U \bar{U}
\eeaa
where
\beaa
U=\mbox{det}(M_n) =\sum_{\sigma}\mbox{sign}(\sigma)\prod_{j=1}^nr_j^{\sigma(j)-1}e^{i(\sigma(j)-1)\theta_j}
\eeaa
and $\sigma$ is a permutation running over all elements in the symmetric group $S_n$.

Note that $\int_{[0,2\pi)}e^{ij\theta}d\theta=0$ for any integer $j\ne 0$. Therefore, any two terms in the sum are orthogonal to each other, that is, for any $\sigma \ne \sigma'$,
\beaa
\int_{[0, 2\pi]^n}\prod_{j=1}^ne^{i(\sigma(j)-1)\theta_j}\cdot \prod_{j=1}^ne^{-i(\sigma'(j)-1)\theta_j}\,d\theta_1\cdots\,d\theta_n=0.
\eeaa
Thus,
\beaa
\int_{[0, 2\pi]^n}U\bar{U}\,d\theta_1\cdots\,d\theta_n=(2\pi)^n\sum_{\sigma}\prod_{j=1}^nr_j^{2(\sigma(j)-1)}.
\eeaa
By integrating out all $\theta_j$'s in (\ref{mean}), we get that the probability density function of  $(R_1, \cdots, R_n)$ is equal to
\beaa
C\cdot (2\pi)^n(r_1\cdots r_n)\sum_{\sigma\in S_n}\prod_{j=1}^nr_j^{2(\sigma(j)-1)}\varphi(r_j)=C\cdot (2\pi)^n\sum_{\sigma\in S_n}\prod_{j=1}^nr_j^{2\sigma(j)-1}\varphi(r_j)
\eeaa
for $r_1\geq 0, \cdots, r_n\geq 0$ and the density is $0$, otherwise. It follows that
\beaa
1 &= & C(2\pi)^n\int_{[0,\infty)^n}\sum_{\sigma\in S_n}\prod_{j=1}^nr_j^{2\sigma(j)-1}\varphi(r_j)\,dr_1\cdots\,dr_n\\
& = & C(2\pi)^n\sum_{\sigma\in S_n}\int_{[0,\infty)^n}\prod_{j=1}^nr_j^{2\sigma(j)-1}\varphi(r_j)\,dr_1\cdots\,dr_n.
\eeaa
For each $\sigma\in S_n$, it is easy to see that the integral is equal to $\prod_{j=1}^n\int_0^{\infty}x^{2j-1}\varphi(x)\,dx.$  We then get the value of $C.$

 Second, the first step says that the density function of $(Z_1, \cdots, Z_n)$ is
\bea\lbl{rain}
f(z_1, \cdots, z_n)=\frac{1}{n!}\cdot \frac{1}{c_0\cdots c_{n-1}}\prod_{1\leq j < k
\leq n}|z_j-z_k|^2\cdot \prod_{j=1}^n\varphi(|z_j|).
\eea
Now write
\beaa
\frac{1}{(c_0\cdots c_{n-1})^{1/2}}\prod_{1\leq j < k
\leq n}(z_j-z_k)=\mbox{det}\begin{pmatrix}
p_0(z_1) & p_0(z_2) & \cdots & p_0(z_n)\\
p_1(z_1) & p_1(z_2) & \cdots & p_1(z_n)\\
\vdots\\
p_{n-1}(z_1) & p_{n-1}(z_2) & \cdots & p_{n-1}(z_n)
\end{pmatrix}
,
\eeaa
where $p_l(z)=z^l/\sqrt{c_{l}}$ for $l=0,1,2\cdots$.  Let $\bd{A}$ be the above matrix. Then,
\bea
  \frac{1}{c_0\cdots c_{n-1}}\prod_{1\leq j < k
\leq n}|z_j-z_k|^2
&=&  \mbox{det}(\bd{A}^*\bd{A}) \nonumber\\
&=&\big(\sum_{k=0}^{n-1}p_k(z_i)\overline{p_k(z_j)}\,\big)_{n\times n}=\mbox{det}((K(z_i, z_j)_{n\times n}).\lbl{sky}
\eea
By the polar transformation,
$\int_{\mathbb{C}^2}z^{j}\bar{z}^k\varphi(|z|)\,dz
= 2\pi\int_0^{\infty}r^{2j+1}\varphi(r)\,dr=c_j$ for $j=k$, and the integral is equal to $0$ for any non-negative integers $j \ne k.$ Hence $\{p_0(z), p_1(z), p_2(z), \cdots\}$ are orthonormal with respect to the measure $\varphi(|z|)\,dz$. By using Exercise 4.1.1 from Hough {\it et al}. (2009), we get the desired conclusion from (\ref{rain}) and (\ref{sky}).  \hfill$\blacksquare$

\medskip

To prove Proposition \ref{Aha_you}, we need some basic facts regarding point processes.  Let $Z_1, \cdots, Z_N$ be random variables with symmetric density $f(z_1,\cdots, z_N)$ with respect to reference measure $\mu$ on $\mathbb{C}$. View them as a point process $\sum_{i=1}^N\delta_{Z_i}$. Then this process has $n$-point correlation function
\beaa
\rho_n(z_1, \cdots, z_n)=\frac{N!}{(N-n)!}\int_{\mathbb{C}^{N-n}}f(z_1,\cdots, z_N)\,\mu(dz_{n+1})\cdots \mu(dz_N)
\eeaa
for $1\leq n <N$ and $\rho_N(z_1, \cdots, z_N)=N!f(z_1,\cdots, z_N)$.
Let $f_n(z_1, \cdots, z_n)$ be the joint density of $Z_1, \cdots, Z_n$. Then,
\begin{equation}\label{f2rho}
f_n(z_1, \cdots, z_n)=\frac{(N-n)!}{N!}\rho_n(z_1, \cdots, z_n)
\end{equation}
for $1\leq n \leq N.$ This means that, for any measurable function $h(z_1, \cdots, z_n)$ with $1\leq n\leq N$, we have
\bea\lbl{face_catch}
Eh(Z_1, \cdots, Z_n)=\frac{(N-n)!}{N!}\int_{\mathbb{C}^n}h(z_1, \cdots, z_n)\rho_n(z_1, \cdots, z_n)\prod_{j=1}^n\mu(dz_j)
\eea
provided
\beaa
\int_{\mathbb{C}^n}|h(z_1, \cdots, z_n)|\,\rho_n(z_1, \cdots, z_n)\prod_{j=1}^n\mu(dz_j) < \infty.
\eeaa
See further details from, for example, Johansson (2005) and (1.2.9) from Hough {\it et al.} (2009). Let $Z_1, \cdots, Z_n$  be  as in Proposition \ref{Laplace}. Then
\bea\lbl{determinantal}
 \rho_k(z_1, \cdots, z_k)=\mbox{det}\, \big(K(z_i, z_j)\big)_{1\leq i, j \leq k} \mbox{ where } K(z, w)=\sum_{k=0}^{n-1}\frac{1}{c_k}(z\bar{w})^k.
\eea
Evidently, $\overline{K(z, w)}=K(w, z)$ for all $(z, w)\in \mathbb{C}^2$. Further, the product of the diagonal entries of $(K(z_i, z_j))_{1\leq i, j \leq k}$ is equal to $\prod^k_{j=1}\rho_1(z_j)$.

\medskip

\noindent\textbf{Proof of Proposition \ref{Aha_you}.} For convenience we now switch the notation $n$ to $N$. So we need to prove
\begin{equation}\label{4thmoment}
\E\Big[\sum^N_{j=1}(h(Z_j)-\E h(Z_j))\Big]^4\le K N^2
\end{equation}
for all $N\ge 1$, where $K$ is a constant not depending on $N$, $\varphi(z)$ or $h(z).$

Obviously, $h(z)=h^{+}(z)-h^{-}(z)$ where $h^{+}(z):=\max\{h(z), 0\}$ and $h^{-}(z):=\min\{h(z), 0\}$. With the trivial bound $(a+b)^4 \leq 8(a^4+b^4)$ for  $(a,b)\in \mathbb{R}^2$, to prove the proposition, we can simply assume $0\leq h(z)\leq 1$ for all $z\in \mathbb{C}.$

Set $e_h=\E h(Z_1)$. Then $e_h \in [0, 1].$ Without loss of generality, we assume $\mu$ is a probability measure, i.e., $\mu(\mathbb{C})=1.$ Review the identity that
\beaa
& & (a_1+\cdots +a_N)^4\\
&=& \sum_{i=1}^Na_i^4 + 6\sum_{i<j}a_i^2a_j^2 +\sum_{i\ne j}a_i^3a_j+  12\sum_{i\ne j \ne k}a_i^2a_ja_k+24\sum_{i<j<k<l}a_ia_ja_ka_l
\eeaa
where $i,j,k,l$ are all between $1$ to $N$.  Set $a_j=h(Z_j)-e_h$,  $1\le j\le N$.  Since $h(z)$ is a bounded function, it
is easy to see that,  to complete the proof of \eqref{4thmoment},
 it suffices to show
\begin{equation}\label{estimate3}
 2N(N-1)(N-2) \E[(h(Z_1)-e_h)^2(h(Z_2)-e_h)(h(Z_3)-e_h)]  \le 6N^2
\end{equation}
and
\begin{equation}\label{estimate4}
 N(N-1)(N-2)(N-3)\E\Big[\prod_{j=1}^4(h(Z_j)-e_h)\Big]
\le 82N^2+18N
\end{equation}
for all $N\geq 4.$
Define $u_k(z_1,\cdots,z_k)= \prod^k_{j=1}\rho_1(z_j)-\rho_k(z_1,\cdots,z_k)$ for any $k\ge 1$.
By using the Hadamard inequality, $\rho_k(z_1,\cdots,z_k)\le \prod^k_{j=1}\rho_1(z_j)$, we know
$u_k(z_1,\cdots,z_k)\ge 0$. Then using \eqref{f2rho}
\[
u_k(z_1,\cdots,z_k)=\Big[N^k\prod^k_{j=1}f_1(z_j)\Big]-N(N-1)\cdots(N-k+1)f_k(z_1,\cdots,z_k).
\]
Since $1-\prod_{j=1}^k(1-x_j)\leq \sum_{j=1}^kx_j$ for all $x_j\in [0,1]$, we obtain
\begin{equation}\label{bound}
\int_{\mathbb{C}^k}u_k(z_1,\cdots,z_k)\prod^k_{j=1}\mu(dz_j)=N^k-N(N-1)\cdots (N-k+1)\le \frac{k(k-1)}{2}N^{k-1}.
\end{equation}
  It follows from
(\ref{f2rho}) that
\beaa
& & 2N(N-1)(N-2) \E[(h(Z_1)-e_h)^2(h(Z_2)-e_h)(h(Z_3)-e_h)]\\
& = & 2\int_{\mathbb{C}^3}(h(z_1)-e_h)^2(h(z_2)-e_h)(h(z_3)-e_h)\rho_3(z_1, z_2, z_3)\prod_{j=1}^3\mu(dz_j)\\
& = & 2\int_{\mathbb{C}^3}(h(z_1)-e_h)^2(h(z_2)-e_h)(h(z_3)-e_h)\rho_1(z_1)\rho_1(z_2)\rho_1(z_3)\prod_{j=1}^3\mu(dz_j)\\
&&-2\int_{\mathbb{C}^3}(h(z_1)-e_h)^2(h(z_2)-e_h)(h(z_3)-e_h)u_3(z_1,z_2, z_3)\prod_{j=1}^3\mu(dz_j).
\eeaa
Since $\rho_1(z_1)\rho_1(z_2)\rho_1(z_2)=N^3f_1(z_1)f_1(z_2)f_1(z_3)$ from \eqref{f2rho}, we have
\[
2\int_{\mathbb{C}^3}(h(z_1)-e_h)^2(h(z_2)-e_h)(h(z_3)-e_h)\rho_1(z_1)\rho_1(z_2)\rho_1(z_3)\prod_{j=1}^3\mu(dz_j)=0.
\]
Thus, we get  from \eqref{bound} and the assumption $h(z)\in[0,1]$ for all $z\in \mathbb{C}$ that
\beaa
& & \Big|2N(N-1)(N-2) \E[(h(Z_1)-e_h)^2(h(Z_2)-e_h)(h(Z_3)-e_h)]\Big|\\
&=&\Big|-2\int_{\mathbb{C}^3}(h(z_1)-e_h)^2(h(z_2)-e_h)(h(z_3)-e_h)u(z_1,z_2,z_3)\prod_{j=1}^3\mu(dz_j)\Big|\\
&\le &2\int_{\mathbb{C}^3}u(z_1,z_2,z_3)\prod_{j=1}^3\mu(dz_j)\\
&=&6N^2,
\eeaa
proving \eqref{estimate3}. Now, we start to prove (\ref{estimate4}).

Define $b_k=\E[\,\prod^k_{j=1}h(Z_j)]$ for $1\le k\le 4$.
Note that
\beaa
\prod_{j=1}^4(h(Z_j)-e_h)&=&e_h^4-e_h^3\sum_{1\le i\le 4}h(Z_i)+e_h^2\sum_{1\le i< j\le 4}h(Z_i)h(Z_j)\\
&&~~~
-e_h\sum_{1\le i< j<k\le 4}h(Z_i)h(Z_j)h(Z_k)+h(Z_1)h(Z_2)h(Z_3)h(Z_4).
\eeaa
Taking expectations on both sides and noting that $b_1=e_h$ we have
\begin{equation}\label{best}
\E\Big[\prod_{j=1}^4(h(Z_j)-e_h)\Big]=-3e_h^4+6e_h^2b_2-4e_hb_3+b_4.
\end{equation}
It follows from
(\ref{face_catch}) that for $2\le k\le 4$
\beaa
b_k&=&\frac{(N-k)!}{N!}\int_{\mathbb{C}^k}\Big[\prod_{j=1}^kh(z_j)\Big]\rho_k(z_1, \cdots, z_k)\prod_{j=1}^k\mu(dz_j)\\
&=&\frac{(N-k)!}{N!}\int_{\mathbb{C}^k}\Big[\prod_{j=1}^kh(z_j)\rho_1(z_j)\Big]\prod_{j=1}^k\mu(dz_j)\\
&&+\frac{(N-k)!}{N!}\int_{\mathbb{C}^k}\Big[\prod_{j=1}^kh(z_j)\Big]\Big(\rho_k(z_1, \cdots, z_k)-\prod_{j=1}^k\rho_1(z_j) \Big)\prod_{j=1}^k\mu(dz_j)\\
&=&\frac{(N-k)!N^k}{N!}e_h^k-
\frac{(N-k)!}{N!}\int_{\mathbb{C}^k}\Big[\prod_{j=1}^kh(z_j)\Big]u_k(z_1, \cdots, z_k)\prod_{j=1}^k\mu(dz_j).
\eeaa
We next evaluate the last integral  for $k=2, 3, 4$.  Denote
\bea\lbl{monkey_bar}
\alpha_k=\int_{\mathbb{C}^k}\Big[\prod_{j=1}^kh(z_j)\Big]u_k(z_1, \cdots, z_k)\prod_{j=1}^k\mu(dz_j).
\eea
Then we have
\begin{equation}\label{estimateofbk}
b_k=\frac{(N-k)!N^k}{N!}e_h^k-
\frac{(N-k)!}{N!}\alpha_k
\end{equation}
for $k=2,3,4$. Obviously, $\alpha_k\ge 0$ for $2\le k\le 4$ since $u_k(z_1, \cdots, z_k)\geq 0$  for all $(z_1, \cdots, z_k) \in \mathbb{C}^k.$ In order to prove (\ref{estimate4}), we need to study $b_2, b_3, b_4$ in (\ref{best}). Based on (\ref{estimateofbk}), it suffices to work on $\alpha_2, \alpha_3, \alpha_4$. We will do so step by step in the following.

\medskip
\noindent{\it Estimate of  $\alpha_2$}. It follows from \eqref{bound} that
\begin{equation}\lbl{commencement}
0\leq \alpha_2\le N.
\end{equation}
Meanwhile, by the determinantal formula (\ref{determinantal}),
since that
\beaa
\rho_2(z_1,z_2)&=&K(z_1,z_1)K(z_2,z_2)-K(z_1,z_2)K(z_2,z_1)\\
&=&\rho_1(z_1)\rho_1(z_2)-|K(z_1,z_2)|^2,
\eeaa
we have
\begin{equation}\label{alpha2}
\alpha_2=\int_{\mathbb{C}^2}h(z_1)h(z_2)|K(z_1,z_2)|^2 \mu(dz_1)\mu(dz_2),
\end{equation}
which will be used later.

\medskip

\noindent{\it Estimate of $\alpha_3$.} Set
\bea\lbl{crossing}
\beta_3:=\int_{\mathbb{C}^3}\Big[\prod^3_{j=1}h(z_j)\Big]K(z_1,z_2)K(z_2,z_3)K(z_3,z_1)
\prod^3_{j=1}\mu(dz_j).
\eea
In this step, we will show
\beaa
\max\{\alpha_3, |\beta_3|\} \leq 3N^2.
\eeaa
It is easily seen from \eqref{bound} that $\alpha_3\le 3N^2$. We now estimate $\beta_3$.  By (\ref{determinantal}) again,
\bea
&&u_3(z_1,z_2, z_3)\nonumber\\
&=&\rho_1(z_1)\rho_1(z_2)\rho_1(z_3)-\rho_3(z_1,z_2,z_3)\nonumber\\
&=&K(z_1,z_1)K(z_2,z_3)K(z_3,z_2)+K(z_1,z_3)K(z_2,z_2)K(z_3,z_1)+K(z_1,z_2)K(z_2,z_1)K(z_3,z_3)\nonumber\\
&&-K(z_1,z_2)K(z_2,z_3)K(z_3,z_1)-K(z_1,z_3)K(z_2,z_1)K(z_3,z_2)\nonumber\\
&=&\rho_1(z_1)|K(z_2,z_3)|^2+\rho_1(z_2)|K(z_1,z_3)|^2
+\rho_1(z_3)|K(z_1,z_2)|^2\label{positive}\\
&&-K(z_1,z_2)K(z_2,z_3)K(z_3,z_1)-K(z_1,z_3)K(z_2,z_1)K(z_3,z_2).\nonumber
\eea
All three functions in \eqref{positive} are nonnegative. For the first term in \eqref{positive} we have from \eqref{alpha2} and then \eqref{f2rho} that
\beaa
&&\int_{\mathbb{C}^3}\Big[\prod^3_{j=1}h(z_j)\Big]\rho_1(z_1)|K(z_2,z_3)|^2\prod^3_{j=1}\mu(dz_j)\\
&=&\int_{\mathbb{C}}h(z_1)\rho_1(z_1)\mu(dz_1)\int_{\mathbb{C}^2}h(z_2)h(z_3)|K(z_2,z_3)|^2\mu(dz_2)\mu(dz_3)\\
&=&\alpha_2\int_{\mathbb{C}}h(z_1)\rho_1(z_1)\mu(dz_1)\\
&=&N\alpha_2e_h.
\eeaa
The same is true for other two terms in (\ref{positive}).  Trivially,
\beaa
\beta_3
=\int_{\mathbb{C}^3}\Big[\prod^3_{j=1}h(z_j)\Big] K(z_1,z_3)K(z_2,z_1)K(z_3,z_2)\prod^3_{j=1}\mu(dz_j).
\eeaa
Therefore, we obtain
\begin{equation}\label{alpha3}
\alpha_3=\int_{\mathbb{C}^3}\Big[\prod^3_{j=1}h(z_j)\Big]u_3(z_1,z_2,z_3)\prod^3_{j=1}\mu(dz_j)=3N\alpha_2e_h-2\beta_3,
\end{equation}
which together with the facts $\alpha_2\leq N$ and $\alpha_3\le 3N^2$ implies
\[
|\beta_3|\le \frac{\alpha_3}{2}+\frac{3N\alpha_2}{2}\le 3N^2.
\]

\medskip

\noindent{\it Estimate of  $\alpha_4$.} This step is a bit involved. The sketch of the proof is as follows.
Since $\rho_4(z_1,z_2,z_3,z_4)$ is the determinant of
$\big(K(z_i,z_j)\big)_{1\le i, j\le 4}$, it can be written as the sum of $24$ terms:
\[
\rho_4(z_1,z_2,z_3,z_4)=\sum_{\sigma}\text{sign}(\sigma)K(z_1,z_{\sigma(1)})K(z_2, z_{\sigma(2)})K(z_3,z_{\sigma(3)})K(z_4,z_{\sigma(4)}),
\]
where $\sigma=(\sigma(1), \sigma(2), \sigma(3), \sigma(4))$ runs over all $24$ permutations of $(1,2,3,4)$.
Excluding $(1,2,3,4)$, all other 23 permutations can be classified into one of the following 4 sets:
\[
D_1=\{(1,2, 4,3), (1,4,3,2), (1,3,2,4),(4,2,3,1), (3,2,1,4), (2,1,3,4)\},
\]
\[
D_2=\{(2,1, 4,3), (3,4,1,2), (4,3,2,1)\},
\]
\[
D_3=\{(4,1,2,3),(4,3, 1,2 ),  (3, 4, 2, 1), (3, 1, 4, 2),  (2,3,4,1), (2, 4, 1, 3)\}
\]
and
\[
D_4=\{(1, 4,2,3), (1, 3,4,2), (4,2,1,3), (3,2,4,1), (4,1,3,2), (2,4, 3, 1), (3, 1,2, 4), (2,3, 1, 4)\}.
\]

Define
\[
T_k(z_1,z_2,z_3,z_4)=\sum_{\sigma\in D_k}K(z_1,z_{\sigma(1)})K(z_2, z_{\sigma(2)})K(z_3,z_{\sigma(3)})K(z_4,z_{\sigma(4)})
\]
for $k=1,2,3,4$. Then
\beaa
\rho_4(z_1,z_2,z_3,z_4)=K(z_1,z_1)K(z_2,z_2)K(z_3,z_3)K(z_4,z_4)+\sum^4_{k=1}(-1)^kT_k(z_1,z_2,z_3,z_4).
\eeaa
This implies
\begin{equation}\label{u4}
u_4(z_1,z_2,z_3,z_4)=\sum^4_{k=1}(-1)^{k+1}T_k(z_1,z_2,z_3,z_4).
\end{equation}
Within each class $D_k$,  all $K(z_1,z_{\sigma(1)})K(z_2, z_{\sigma(2)})K(z_3,z_{\sigma(3)})K(z_4,z_{\sigma(4)})$
contribute equally to the integral $\int_{\mathbb{C}^4}\big[\prod^4_{j=1}h(z_j)\big]T_k(z_1,z_2,z_3,z_4)\prod^4_{j=1}\mu(dz_j)$.
We have
\beaa
&&\int_{\mathbb{C}^4}\big[\prod^4_{j=1}h(z_j)\big]T_1(z_1,z_2,z_3,z_4)\prod^4_{j=1}\mu(dz_j)\\
&=&6\int_{\mathbb{C}^4}\big[\prod^4_{j=1}h(z_j)\big]K(z_1,z_1)K(z_2,z_2)K(z_3,z_4)K(z_4,z_3)\prod^4_{j=1}\mu(dz_j)\\
&=&6\int_{\mathbb{C}^4}\big[\prod^4_{j=1}h(z_j)\big]\rho_1(z_1)\rho_1(z_2)K(z_3,z_4)K(z_4,z_3)\prod^4_{j=1}\mu(dz_j)\\
&=&6N^2e_h^2\alpha_2
\eeaa
by using (\ref{alpha2})
and the definition that $e_h=\E h(Z_1)$.  For $T_2$, we have
\beaa
&&\int_{\mathbb{C}^4}\big[\prod^4_{j=1}h(z_j)\big]T_2(z_1,z_2,z_3,z_4)\prod^4_{j=1}\mu(dz_j)\\
&=&3\int_{\mathbb{C}^4}\big[\prod^4_{j=1}h(z_j)\big]K(z_1,z_2)K(z_2,z_1)K(z_3,z_4)K(z_4,z_3)\prod^4_{j=1}\mu(dz_j)\\
&=&3\int_{\mathbb{C}^2}h(z_1)h(z_2)|K(z_1,z_2)|^2\mu(dz_1)\mu(dz_2)\\
&&~~~~\times \int_{\mathbb{C}^2}h(z_3)h(z_4)|K(z_3,z_4)|^2\mu(dz_3)\mu(dz_4)\\
&=&3\alpha_2^2\\
&\le &3N^2
\eeaa
by (\ref{commencement}). Noting that $\prod^4_{j=1}\mu(dz_j)$ is a probability measure, we have from the Cauchy-Schwarz inequality and the fact $0\leq h(z)\leq 1$ for all $z\in \mathbb{C}$ that
\beaa
&&\int_{\mathbb{C}^4}\big[\prod^4_{j=1}h(z_j)\big]T_3(z_1,z_2,z_3,z_4)\prod^4_{j=1}\mu(dz_j)\\
&=&6\int_{\mathbb{C}^4}\big[\prod^4_{j=1}h(z_j)\big]K(z_1,z_4)K(z_2,z_1)K(z_3,z_2)K(z_4,z_3)\prod^4_{j=1}\mu(dz_j)\\
&\leq&6\left(\int_{\mathbb{C}^4}\Big[\prod^4_{j=1}h(z_j)\Big]|K(z_1,z_4)K(z_3,z_2)|^2\prod^4_{j=1}\mu(dz_j)\right)^{1/2}\\
&&~~~\times\left(\int_{\mathbb{C}^4}\Big[\prod^4_{j=1}h(z_j)\Big]|K(z_1,z_4)K(z_3,z_2)|^2\prod^4_{j=1}\mu(dz_j)\right)^{1/2}\\
&=&6\int_{\mathbb{C}^4}\Big[\prod^4_{j=1}h(z_j)\Big]|K(z_1,z_4)|^2|K(z_3,z_2)|^2\prod^4_{j=1}\mu(dz_j)\\
&=&6\Big[\int_{\mathbb{C}^2}h(z_1)h(z_4)|K(z_1,z_4)|^2\mu(dz_1)\mu(dz_4)\Big]^2\\
&=&6\alpha_2^2\\
&\le &6N^2
\eeaa
by \eqref{alpha2}  and (\ref{commencement}).

We next estimate the term on $T_4$. In fact,
\beaa
&&\int_{\mathbb{C}^4}\big[\prod^4_{j=1}h(z_j)\big]T_4(z_1,z_2,z_3,z_4)\prod^4_{j=1}\mu(dz_j)\\
&=&8\int_{\mathbb{C}^4}\big[\prod^4_{j=1}h(z_j)\big]K(z_1,z_1)K(z_2,z_4)K(z_3,z_2)K(z_4,z_3)\prod^4_{j=1}\mu(dz_j)\\
&=&8\int_{\mathbb{C}^4}\big[\prod^4_{j=1}h(z_j)\big]\rho_1(z_1)K(z_2,z_4)K(z_3,z_2)K(z_4,z_3)\prod^4_{j=1}\mu(dz_j)\\
&=&8Ne_h\int_{\mathbb{C}^3}\big[\prod^4_{j=2}h(z_j)\big]K(z_2,z_4)K(z_3,z_2)K(z_4,z_3)\prod^4_{j=2}\mu(dz_j)\\
&=&8Ne_h\beta_3,
\eeaa
by (\ref{crossing}).

Now multiplying $\prod^4_{j=1}h(z_j)$ on both sides of \eqref{u4} and integrating with respect to the measure $\prod^4_{j=1}\mu(dz_j)$
we obtain from (\ref{monkey_bar}) that
\begin{equation}\label{alpha4}
\alpha_4=6N^2e_h^2\alpha_2 - 8Ne_h\beta_3+d_{23},
\end{equation}
where
$|d_{23}|\le 9N^2$. From  \eqref{estimateofbk}, \eqref{alpha3} and \eqref{alpha4} we see that
\beaa
& & b_2=\frac{N e_h^2}{N-1}-\frac{\alpha_2,}{N(N-1)};\\
& & b_3=\frac{N^2 e_h^3}{(N-1)(N-2)}-\frac{3Ne_h\alpha_2-2\beta_3}{N(N-1)(N-2)};\\
& & b_4=\frac{N^3e_h^4}{(N-1)(N-2)(N-3)}-\frac{6N^2e_h^2\alpha_2-8Ne_h\beta_3+d_{23}}{N(N-1)(N-2)(N-3)}.
\eeaa
Then it follows from \eqref{best} that
\beaa
&&\E\Big[\prod_{j=1}^4(h(Z_j)-e_h)\Big]\\
&=&-3e_h^4+6e_h^2b_2-4e_hb_3+b_4\\
&=&\Big(-3+\frac{6N}{N-1}-\frac{4N^2}{(N-1)(N-2)}+\frac{N^3}{(N-1)(N-2)(N-3)}\Big)e_h^4\\
&&+\Big(-\frac{6}{N(N-1)}+\frac{12N}{N(N-1)(N-2)}-\frac{6N^2}{N(N-1)(N-2)(N-3)}\Big)e_h^2\alpha_2\\
&&-\Big(\frac{8}{N(N-1)(N-2)} -\frac{8N}{N(N-1)(N-2)(N-3)}\Big)e_h\beta_3\\
&&-\frac{d_{23}}{N(N-1)(N-2)(N-3)}\\
&=&\frac{(3N+18)e_h^4}{(N-1)(N-2)(N-3)}-\frac{(6N+36)e_h^{2} \alpha_2}{N(N-1)(N-2)(N-3)}\\
&&+\frac{24 e_h\beta_3}{N(N-1)(N-2)(N-3)}-\frac{d_{23}}{N(N-1)(N-2)(N-3)}\\
&\le&\frac{90N^2+54N}{N(N-1)(N-2)(N-3)}
\eeaa
for $N\geq 4$ by  the facts $0\leq e_h\leq 1$, $0\leq \alpha_2\leq N$, $|\beta_3|\leq 3N^2$ and $|d_{23}| \leq 9N^2$.  This proves \eqref{estimate4}. The proof is then completed. \hfill$\blacksquare$

\medskip

\noindent{\it Proof of Proposition \ref{lem1}.} We will only need to prove (ii). In fact, conclusion (i)  follows from (ii)
since $h(z)=h(|z|e^{i\theta})$ with $z=|z|e^{i\theta}$. By Proposition \ref{Laplace} and (\ref{f2rho}), the density function  $f_1(z_1)$ of $Z_1$ is given by
\[
f_1(z)\varphi(|z|)=\frac{1}{n}K(z,\bar{z})\varphi(|z|)=\frac1n\sum^{n-1}_{k=0}\frac{|z|^{2k}\varphi(|z|)}{c_k}.
\]
Write $z=x+yi=re^{i\theta}$ with $r\geq 0$ and $\theta\in [0, 2\pi)$. Then $(x, y)=(r\cos\theta, r\sin\theta)$. The Jacobian for the transformation is known to be $r$.
By (\ref{densityofYj}), (\ref{jilin}) and Proposition \ref{Laplace},
this implies that the joint density function of $|Z_1|$ and $\Theta_1$ is given by
\[
f(r, \theta)=r\cdot \frac1n\sum^{n-1}_{k=0}\frac{r^{2k}\varphi(r)}{c_k}=\frac1n\sum^{n-1}_{k=0}\frac{r^{2k+1}\varphi(r)}{c_k}=\frac{1}{2\pi}P_n(r)
\]
for $r\ge 0$ and $\theta\in [0,2\pi)$. Therefore, $|Z_1|$ and $\Theta_1$ are independent,  the density function of
$|Z_1|$ is $P_n(r)$, and $\Theta_1$ is uniformly distributed over $[0,2\pi)$. The conclusion (\ref{Coke_pepsi}) follows immediately.
\hfill$\blacksquare$

\medskip

\noindent\textbf{Proof of Lemma \ref{Catroll}}. Let $Y_{nj}$,  $1\le j\le n$, be independent random variables such that
$Y_{nj}$ has a density function proportional to $y^{2j-1}w_m(y)$ for each $j$, {\it where $w_m(\cdot)$ is as in (\ref{up_stairs})}. Then it follows from Lemma~\ref{independence} that $g(|Z_1|^2, \cdots, |Z_n|^2)$  and $g(Y_{n1}^2, \cdots, Y_{nn}^2)$ are identically distributed. In the proof of their Lemma 2.4, Jiang and Qi (2017) show that $Y_{nj}^2$ has the same distribution as that of $\prod_{r=1}^ms_{j,r}$ for each $1\le j\le n$. This yields
the  desired conclusion.     \hfill$\blacksquare$

\medskip

\noindent\textbf{Proof of Lemma \ref{new_loving}}. Define $\rho^{(l)}(s)=(1-s)^{l-1}I(0\le s<1)$ for $l\geq 1$. Then we have $w_1^{(l)}(s)=(l/\pi)\rho^{(l)}(s^2)$.
Set $\rho_1^{(l_1)}(s)=\rho^{(l_1)}(s)$ and define $\rho_m^{(l_1, \cdots, l_m)}(s)$ recursively by
\begin{equation}\label{holloween}
\rho_m^{(l_1, \cdots, l_m)}(s)=\int^1_0\rho^{(l_m)}(u)\rho_{m-1}^{(l_1, \cdots, l_{m-1})}(\frac{s}u)\frac{du}{u}
\end{equation}
for $m\geq 2$ and positive integers $l_1, \cdots, l_m$.
Evidently the support of $\rho_m^{(l_1, \cdots, l_m)}(s)$ is $[0,1]$. By induction, it is easy to verify from (\ref{salt_sugar}) that
\begin{equation}\label{w2rho}
w_m^{(l_1, \cdots, l_m)}(s)=\frac{1}{\pi}(\prod^m_{j=1}l_j)\rho_m^{(l_1, \cdots, l_m)}(s^2)
\end{equation}
for $m\ge 1$.
Define
\[
\theta^{(l_1, \cdots, l_m)}_m(t)=\int^1_0s^{t-1}\rho_m^{(l_1, \cdots, l_m)}(s)ds,\ \ t>0.
\]
Obviously,  $\theta^{(l_1)}_1(t)=B(t, l_1)$. Now, for any $m>1$, we have from \eqref{holloween} that
\beaa
\theta^{(l_1, \cdots, l_m)}_m(t)
&=&\int^1_0\int^1_0s^{t-1}\rho^{(l_m)}(u)\rho_{m-1}^{(l_1, \cdots, l_{m-1})}(\frac{s}u)\frac{duds}{u}\\
&=&\int^1_0u^{t-1}\rho^{(l_m)}(u)\Big(\int^1_0(\frac{s}{u})^{t-1}\rho_{m-1}^{(l_1, \cdots, l_{m-1})}(\frac{s}u)ds\Big)\frac{du}u\\
&=&\int^1_0u^{t-1}\rho^{(l_m)}(u)\Big(\int^{1/u}_0y^{t-1}\rho_{m-1}^{(l_1, \cdots, l_{m-1})}(y)dy\Big)du.
\eeaa
Keeping in mind that the support of $\rho_m^{(l_1, \cdots, l_m)}(s)$ is $[0,1]$, the above is identical to
\beaa
&&\int^1_0u^{t-1}\rho^{(l_m)}(u)\Big(\int^{1}_0y^{t-1}\rho_{m-1}^{(l_1, \cdots, l_{m-1})}(y)dy\Big)du\\
&=&\theta^{(l_1, \cdots, l_{m-1})}_{m-1}(t) \int^1_0u^{t-1}(1-u)^{l_m-1}du\\
&=&\theta^{(l_1, \cdots, l_{m-1})}_{m-1}(t)B(t, l_m).
\eeaa
We thus conclude from the recursive formula that for any $m\ge 1$,
\begin{equation}\label{candy}
\theta^{(l_1, \cdots, l_m)}_m(t)=\int^1_0s^{t-1}\rho_m^{(l_1, \cdots, l_m)}(s)ds=\prod^m_{r=1}B(t, l_r).
\end{equation}


Let $q_j(y)$ be a density function proportional to $y^{j-1}w_{m}^{(l_1, \cdots, l_m)}(y^{1/2})$ for $1\leq j \leq n$. Thus, $q_j(y)$ is also proportional to $y^{j-1}\rho_{m}^{(l_1, \cdots, l_m)}(y)$ from \eqref{w2rho}, that is, for some $c_j>0$, $q_j(y)=c_j^{-1} y^{j-1}\rho_{m}^{(l_1, \cdots, l_m)}(y)$. Let $Y_{n,j}$ be a random variable such that the density function of $Y_{n,j}^2$ is $q_j(y)$.
%
Since $q_j(y)$ is a  density, then $c_j=\int^1_0y^{j-1}\rho_{m}^{(l_1, \cdots, l_m)}(y)dy = \theta_m^{(l_1, \cdots, l_m)}(j)=\prod^m_{r=1}B(j, l_r)$ from \eqref{candy}, and hence
\[
q_j(y)=\frac{y^{j-1}\rho_{m}^{(l_1, \cdots, l_m)}(y)}{\prod^m_{r=1}B(j, l_r)}, \ \ 0<y<1.
\]
Denote by $M_j(t)$ the moment generating function of $\log Y_{n,j}^2$. Then,
\begin{eqnarray}\label{logY}
M_j(t)=\E(e^{t\log Y_{n,j}^2}) &= & \int^1_0y^tq_j(y)dy \nonumber\\
&=&\frac{\int^1_0y^{j+t-1}\rho_{m}^{(l_1, \cdots, l_m)}(y)dy}
{\prod^m_{r=1}B(j, l_r)}\nonumber\\
&=&\prod^m_{r=1}\frac{B(j+t, l_r)}{B(j, l_r)}
\end{eqnarray}
for $t>-j$, which is the same as $\prod^m_{r=1}E\exp\{t\log s_{j,r}\}$, where $\{s_{j, r},\,1\leq j\leq n,\, 1\le r\le m\}$ are independent random variables and $s_{j, r}$ has the Beta  density $\frac{1}{B(j,l_r)}y^{j-1}(1-y)^{l_r-1}I(0\leq y\leq 1)$ for each $j$ and $r$.
Then, $Y_{n,j}^2$ and $\prod^m_{r=1}s_{j, r}$ have the same distribution.
This and Lemma \ref{independence} lead to the desired conclusion.
\hfill$\blacksquare$

\subsection{Proof of Theorem~\ref{nonlinear}}\lbl{yes_right}

\noindent{\bf Proof of Theorem~\ref{nonlinear}.} Write $Z_j=R_je^{i\Theta_j}$ with $\Theta_j \in [0, 2\pi)$ for $1\leq j \leq n$. We need to show that for any continuous function $u( \theta, r)$ with $0\le u(\theta, r)\le 1$ for all $\theta\in [0, 2\pi)$ and $r\geq 0$,
\begin{equation}\label{goal}
\frac{1}{n}\sum^n_{k=1}u(\Theta_k, h_n(R_k))\to \frac{1}{2\pi}\int^\infty_0\int^{2\pi}_0u(\theta, r)d\theta \nu(dr)\ \ \ \ a.s.
\end{equation}
as $n\to\infty$.   Obviously, $(\Theta_k, R_k)$ for $1\leq k \leq n$ have the same distribution. First, by the Markov inequality and Proposition~\ref{Aha_you},
\beaa
&&P\Big(\big|\frac{1}{n}\sum^n_{k=1}u(\Theta_k, h_n(R_k))-\E u(\Theta_1, h_n(R_1))\big|\ge \varepsilon\Big)\\
&\le& \frac{\E\big[\sum^n_{k=1}\big(u(\Theta_k, h_n(R_k))-\E u(\Theta_1, h_n(R_1))\big)\big]^4}{n^4\varepsilon^4}\\
&\le &\frac{C}{n^2\varepsilon^4}
\eeaa
for every  $\varepsilon>0$, where $C>0$ is a constant not depending on $n$ or $\epsilon$. This implies that
\[
\sum^\infty_{n=1}P\Big(\big|\frac{1}{n}\sum^n_{k=1}\big(\Theta_k, u(h_n(R_k))-\E u(\Theta_1, h_n(R_1))\big)\big|\ge \varepsilon\Big)<\infty.
\]
We conclude from the Borel-Cantelli lemma that, with probability one,

\begin{equation}\label{goal2}
\frac{1}{n}\sum^n_{k=1}\big( u(\Theta_k, h_n(R_k))-\E u(\Theta_1, h_n(R_1))\big)\to 0
\end{equation}
as $n\to\infty.$

Note that $G(r)=\frac{1}{2\pi}\int^{2\pi}_0u(\theta, r)d\theta$ is bounded and continuous in $r\in [0, \infty).$
Since $ \nu_n$ converges weakly to
$\nu$ with probability one,  we have
\[
\frac{1}n\sum^n_{k=1}G(h_n(Y_k))\to \int^\infty_0G(r) \nu(dr)
\]
with probability one. This implies
\[
\E\Big|\frac{1}n\sum^n_{k=1}G(h_n(Y_k))-\int^\infty_0G(r) \nu(dr)\Big|\to 0
\]
via the bounded convergence theorem.  Hence
\begin{equation*}
\frac{1}n\E\sum^n_{k=1}G(h_n(Y_k))-\int^\infty_0G(r) \nu(dr)\to 0,
\end{equation*}
which together with Proposition~\ref{lem1} yields
\[
\E u(\Theta_1, h_n(R_1))=\int^{\infty}_0G(h_n(r))P_n(r)dr=\frac{1}n\E\sum^n_{k=1}G(h_n(Y_k))\to \int^\infty_0G(r) \nu(dr).
\]
This and \eqref{goal2} imply \eqref{goal}. The proof is then completed.

Now we prove the conclusion for $(\mu_n^*, \nu_n^*, \mu^*, \nu^*)$. It suffices to show that, for any continuous $f(z)$ with $0\le f(z)\le 1$ for every $z\in \mathbb{C}$,
\begin{equation}\label{goal-again}
\frac{1}{n}\sum^n_{k=1}f\big(\frac{Z_k}{a_n}\big)\to \frac{1}{2\pi}\int^\infty_0\int^{2\pi}_0f(re^{i\theta})d\theta\nu^*(dr)
\end{equation}
with probability one.

Define
$g(\theta, r)=f(re^{i\theta})$ and $h_n(r)=r/a_n$.  By Theorem~\ref{nonlinear}, \eqref{goal} holds. It follows that, with probability one,
 \beaa
 \frac{1}{n}\sum^n_{k=1}f\big(\frac{Z_k}{a_n}\big)&=&\frac{1}{n}\sum^n_{k=1}g(\Theta_k, h_n(R_k))\\
 &\to&
\frac{1}{2\pi}\int^\infty_0\int^{2\pi}_0g(\theta, r)d\theta\nu^*(dr)\\
&=&\frac{1}{2\pi}\int^\infty_0\int^{2\pi}_0f(re^{i\theta})d\theta\nu^*(dr),
 \eeaa
completing the proof of \eqref{goal-again}.  \hfill$\blacksquare$

%
%

\subsection{Proof of Theorem \ref{Theorem_product}}

We first need a technical lemma as follows.

\begin{lemma}\label{lem2} Suppose $\{h_n(x);\, n\geq 1\}$ are measurable functions defined on $[0,\infty)$ and $\nu_n$ are defined as in \eqref{mun}. Let $Y_1, \cdots, Y_n$ be as in Lemma~\ref{independence} and $\nu$ be a probability measure on $\mathbb{R}.$
Then  $\nu_n\rightsquigarrow\nu$ if and only if
\begin{equation*}
\lim_{n\to\infty}\frac{1}{n}\sum^n_{j=1}P(h_n(Y_j)\le r)=G(r)
\end{equation*}
for every continuous point $r$ of $G(r)$, where $G(r):=\nu((-\infty, r]),\, r \in \mathbb{R}$.
\end{lemma}

\noindent{\bf Proof.} Let $C_G$ denote the set of all continuity points of $G$.
Note that $\nu_n$ converges weakly to $\nu$ with probability one if and only if $\nu_n((-\infty, r])\to
\nu((-\infty, r])$ with probability one for any $r\in C_{G}$,  that is, for  all $r\in C_G$
\begin{equation}\label{empirical}
\frac{1}{n}\sum^n_{j=1}I(h_n(Y_j)\le r)\to G(r)
\end{equation}
with probability one. Since $Y_1, \cdots, Y_n$ are independent
random variables,
\[
\nu_n((-\infty, r])=\frac{1}{n}\sum^n_{j=1}I(h_n(Y_j)\le r),
\]
which is the average of $n$ independent bounded random variables. By calculating the fourth moment,  applying the
Chebyshev inequality and then the Borel-Cantelli lemma we can show
that  for any $r \in \mathbb{R}$,
\beaa
& &  \frac{1}{n}\sum^n_{j=1}I(h_n(Y_j)\le r)-\frac1n\sum^n_{j=1}P(h_n(Y_j)\le r)\\
&= &
 \frac{1}{n}\sum^n_{j=1}\big[ I(h_n(Y_j)\le r)-P(h_n(Y_j)\le r)\big]\to 0
\eeaa
with probability one. This and (\ref{empirical}) imply the desired conclusion.
  \hfill$\blacksquare$

\medskip

\noindent{\bf Proof of Theorem~\ref{Theorem_product}.} Let $h_n(y)=y^{2/m_n}/n$, $y\ge 0$. By applying Theorem \ref{nonlinear} and Lemma~\ref{lem2} it suffices to show that
\[
\frac{1}{n}\sum^n_{j=1}P(|Z_j|^{2/m_n}/n\le y)\to y, ~~~~y\in (0,1)
\]
  which is equivalent to
\begin{equation}\label{ohyeah}
\frac{1}{n}\sum^n_{j=1}P\Big(\frac{1}{n}\big(\prod_{r=1}^{m_n}s_{j,r}\big)^{1/m_n} \le y\Big)\to y, ~~~~y\in (0,1)
\end{equation}
by Lemma \ref{Catroll}, where $s_{j,r}$'s are as in Lemma \ref{Catroll}. Define $\eta(x)=x-1-\log x$ for $x>0$. Since $\eta(x)=\int^x_1\frac{s-1}{s}ds$, it is easy to verify that
\begin{equation}\label{eta-est}
0\le \eta(x)\le \frac{(x-1)^2}{2\min(x,1)}, ~~~~~~x>0.
\end{equation}
Set $W_j=\prod^{m_n}_{r=1}s_{j, r}$ for $1\leq j\leq n.$ Then,
\[
\log W_j=\sum^{m_n}_{r=1}\log s_{j, r}
\]
for each $j$. By using the expression $\log x=x-1-\eta(x)$ we can rewrite $\log W_{j}$ as
\begin{eqnarray*}
\log W_j&=&\sum^{m_n}_{r=1}\log \frac{s_{j,r}}{j}+m_n\log j\\
&=&\frac1j\sum^{m_n}_{r=1}(s_{j,r}-j)+m_n\log j-\sum^{m_n}_{r=1}\eta\Big(\frac{s_{j,r}}{j}\Big).
\end{eqnarray*}
Write $T_j=\sum^{m_n}_{r=1}s_{j,r}$ for each $1\le j\le n$. Then $T_j$ is the sum of $jm_n$ i.i.d.
random variables with the exponential distribution of mean $1$. Hence $\E(T_j)=jm_n$ and $\Var(T_j)=jm_n$ for $1\le j\le n$. From the above equations we have
\begin{equation}\label{Tsum}
\log W_j=\frac{1}{j}(T_j-jm_n)+m_n\log j-\sum^{m_n}_{r=1}\eta\Big(\frac{s_{j,r}}{j}\Big).
\end{equation}
Since $s_{j, r}$ has the Gamma  density $y^{j-1}e^{-y}I(y>0)/(j-1)!$,
the moment generating functions of $\log s_{j,r}$ is
\begin{equation*}
m_j(t)=\E(e^{t\log s_{j,r}})=\frac{1}{\Gamma(j)}\int^\infty_0y^ty^{j-1}e^{-y}dy=\frac{\Gamma(j+t)}{\Gamma(j)}
\end{equation*}
for $t>-j.$ Therefore,
\bea\lbl{ice_coffee}
\E(\log s_{j,r})=\frac{d}{dt}m_j(t)\Big|_{t=0}=\frac{\Gamma'(j)}{\Gamma(j)}{\red:=}\psi(j).
\eea
The function $\psi$ is the so-called Digamma function in the literature. By Formulas 6.3.18
from Abramowitz and Stegun (1972),
\begin{equation}\label{psi}
\psi(x)=\log x-\frac{1}{2x}+O\Big(\frac{1}{x^2}\Big)~~~~~~\mbox{as  }x\to +\infty.
\end{equation}
Because $\E(\log W_j)=\sum^{m_n}_{r=1}\E\log(s_{j,r})= m_n\psi(j)$,  we see that
\begin{equation}\label{meanoferror}
\sum^{m_n}_{r=1}\E\eta\Big(\frac{s_{j,r}}{j}\Big)=m_n(\log j-\psi(j)).
\end{equation}
Now we fix $y\in (0,1)$.  Write $S_j= \sum^{m_n}_{r=1}\eta(\frac{s_{j,r}}{j})$.
It follows from \eqref{Tsum} that
\begin{eqnarray}\label{key}
&&\frac{1}{n}\sum^n_{j=1}P\Big(\frac{1}{n}\big(\prod_{r=1}^{m_n}s_{j,r}\big)^{1/m_n} \le y\Big)\nonumber\\
&=&\frac{1}{n}\sum^n_{j=1}P\Big(\log W_j\le m_n\log (ny)\Big)\nonumber\\
&=&\frac{1}{n}\sum^n_{j=1}P\Big(\frac{1}{j}(T_j-jm_n)+m_n\log j-S_j\le m_n\log(ny)\Big)\nonumber\\
&=&\frac{1}{n}\sum^n_{j=1}P\Big(\frac{1}{\sqrt{jm_n}}(T_j-jm_n)\le \sqrt{jm_n}\log(\frac{ny}{j})+\sqrt{\frac{j}{m_n}}S_j\Big).
\end{eqnarray}

For any fixed small number $\varepsilon\in (0,1/2)$ such that $y(1+\varepsilon)<1$, define integers $j_n^+=[ny(1+\varepsilon)]+1$ and $j_n^-=[ny/(1+\varepsilon)]$, where $[x]$ denotes the integer part of $x$.   Obviously we have
\begin{equation}\label{part1}
\frac{ny}{j}\ge 1+\varepsilon~~~~~~~\mbox{ if }  j\le j_n^-
\end{equation}
and
\begin{equation}\label{part2}
\frac{ny}{j}\le \frac1{1+\varepsilon}~~~~~~~\mbox{ if } j\ge j_n^+.
\end{equation}
Since $T_j$ is the sum of $jm_n$ i.i.d.
random variables with  both mean and variance equal to $1$. Then $\Var(\frac{1}{\sqrt{jm_n}}(T_j-jm_n))=1.$  From \eqref{part1} and the Chebyshev inequality,
\bea
\frac{1}{n}\sum^{j_n^-}_{j=1}P\Big(\frac{1}{\sqrt{jm_n}}(T_j-jm_n)\ge \sqrt{jm_n}\log\big(\frac{ny}{j}\big)\Big)
&\le &\frac{1}{n}\sum^{j_n^-}_{j=1}\frac{1}{jm_n(\log(1+\varepsilon))^2} \nonumber\\
&=&O\big(\frac{\log n}{nm_n}\big)\to 0 \lbl{Stanford_U}
\eea
as $n\to\infty.$ This implies
\beaa
\lim_{n\to\infty}\frac{1}{n}\sum^{j_n^-}_{j=1}P\Big(\frac{1}{\sqrt{jm_n}}(T_j-jm_n)\le \sqrt{jm_n}\log\big(\frac{ny}{j}\big)\Big) = \frac{y}{1+\epsilon}.
\eeaa
Note that $S_j\ge 0$ from \eqref{eta-est}.  We obtain by  \eqref{key}
\bea
&&\liminf_{n\to\infty}\frac{1}{n}\sum^n_{j=1}P\Big((\prod_{r=1}^{m_n}s_{j,r})^{1/m_n}/n \le y\Big) \nonumber\\
&\ge&\liminf_{n\to\infty}\frac{1}{n}\sum^n_{j=1}P\Big(\frac{1}{\sqrt{jm_n}}(T_j-jm_n)\le \sqrt{jm_n}\log\big(\frac{ny}{j}\big)\Big) \nonumber\\
&=&\frac{y}{1+\varepsilon}. \lbl{sitting_cat}
\eea
From \eqref{part2} we have
\beaa
&&\frac{1}{n}\sum^n_{j=j_n^+}P\Big(\frac{1}{\sqrt{jm_n}}(T_j-jm_n)\le \sqrt{jm_n}\log(\frac{ny}{j})+\sqrt{\frac{j}{m_n}}S_j\Big)\\
&\le &\frac{1}{n}\sum^n_{j=j_n^+}\Big[P\Big(\frac{1}{\sqrt{jm_n}}(T_j-jm_n)\le -\sqrt{jm_n}\log(1+\varepsilon)+\sqrt{\frac{j}{m_n}}S_j, S_j\le \frac{m_n}{2}\log(1+\varepsilon)\Big)\\
&&~~~~~~~~~~~~~~~~~~~~~~~~~~~~~~~~~~~~~~~~~~~~~~~~~~~~~~~~~~~~~~~+P\big( S_j>\frac{m_n}{2}\log(1+\varepsilon)\big)\Big]\\
&\le &\frac{1}{n}\sum^n_{j=j_n^+}P\Big(\frac{1}{\sqrt{jm_n}}(T_j-jm_n)\le -\frac12\sqrt{jm_n}\log(1+\varepsilon)\Big)
     +\frac{1}{n}\sum^n_{j=j_n^+}P( S_j>\frac{m_n}{2}\log(1+\varepsilon)).
\eeaa
By the same argument as in (\ref{Stanford_U}), the first sum above goes to zero. Further, by the Markov inequality, \eqref{meanoferror} and then \eqref{psi}, the last sum is controlled by
\beaa
     \frac{1}{n}\sum^n_{j=j_n^+}\frac{\E(S_j)}{\frac{m_n}{2}\log(1+\varepsilon)}
&\le &
     \frac{1}{n}\sum^n_{j=j_n^+}\frac{m_n(\log j-\psi(j))}{\frac{m_n}{2}\log(1+\varepsilon)}\\
&=&
O(\frac1n)\sum^n_{j=j_n^+}\frac1j \to  0
\eeaa
 as $n\to\infty$ since $\sum^n_{j=j_n^+}\frac1j=O(1)$. 
These and  \eqref{key} imply
 \beaa
&&\limsup_{n\to\infty}\frac{1}{n}\sum^n_{j=1}P\Big((\prod_{r=1}^{m_n}s_{j,r})^{1/m_n}/n \le y\Big)\\
&=& \limsup_{n\to\infty}\frac{1}{n}\sum^{j_n^{+}}_{j=1}P\Big(\frac{1}{\sqrt{jm_n}}(T_j-jm_n)\le \sqrt{jm_n}\log(\frac{ny}{j})+\sqrt{\frac{j}{m_n}}S_j\Big)\\
&\le& \limsup_{n\to\infty}\frac{j_n^{+}}{n}= y(1+\varepsilon).
\eeaa
By taking $\varepsilon\downarrow 0$ to the above and (\ref{sitting_cat}) we get
\eqref{ohyeah}. The proof is completed.
\hfill$\blacksquare$

\subsection{Proof of Theorem~\ref{pie}}\lbl{god_bite}

Let $\{s_{j,r},\,1\le j\le n, 1\leq r\leq m\}$ be independent random variables and $s_{j,r} \sim \,\mbox{Beta} (j, l_r)$, that is, $s_{j,r}$ has the Beta  density $\frac{1}{B(j,l_r)}y^{j-1}(1-y)^{l_r-1}I(0\leq y\leq 1)$ for each $j$ and $r$.
Define
\begin{equation}\label{rep1}
Y_{n,j}^2=\prod^m_{r=1}s_{j,r},  ~~~1\le j\le n.
\end{equation}
Let $[x]$ denote the integer part of $x$ and ``$\overset{p}\to 0$" indicate that ``converges to zero in probability".

We start with an auxiliary result before proving Theorem \ref{pie}.
\begin{lemma}
Let $G_{j}(x)=P(Y_{n,j}^2\le x)$ for  $x\in [0,1]$ and $1\leq j \leq n$. Then,
\begin{equation}\label{stochastic-order}
G_1(x)\ge G_2(x)\ge \cdots\ge G_n(x), \ \ x\in[0,1].
\end{equation}
Further, assume the conditions in Theorem \ref{pie} hold. If (\ref{C_1_condition})  or (\ref{C_2_condition}) is true, then
\begin{equation}\label{stability}
\frac{1}{\gamma_n}\log\frac{Y_{n, [nx]}^2}{b_n}-\log F(x)\overset{p}\to 0,\ \ x\in (0,1),
\end{equation}
as  $n\to\infty$, where $b_n$ and $\gamma_n$ are as in Theorem \ref{pie}.
\end{lemma}

 Note that the two assertions in the above lemma are not directly connected. We put them together simply because they are all about the $Y_{n,j}$'s.

\medskip

\noindent\textbf{Proof}. We first prove \eqref{stochastic-order}.
Let $X_i$ and $Y_i$ be independent positive random variables for $i=1, 2$. If $P(X_1\le x)\ge P(X_2\le x)$ and $P(Y_1\le x)\ge P(Y_2\le x)$ for all $x>0$, then one can easily show that $P(X_1Y_1\le x)\ge P(X_2Y_2\le x)$ for all $x>0$.  Combining this fact
and \eqref{rep1},   it suffices to show that for each $1\le r\le m$ and $x>0$,
$P(s_{j, r}\le x)$ is non-increasing in $j$, where $s_{j, r}$ is as in (\ref{rep1}).

Let $V_j$, $j\ge 1$ be i.i.d. random variables uniformly distributed over $(0,1)$.  For each $n\ge 1$, let
$V_{1:n}\le V_{2:n}\le \cdots\le V_{n:n}$ denote the order statistics of $V_1, \cdots,  V_n$.  It is well-known that
$V_{j:n}$ has a Beta$(j, n-j+1)$ distribution, see e.g., Balakrishnan and Cohen (1991).


It is easy to see that $V_{1:n}\le V_{2: n+1}\le \cdots\le V_{j: n+j-1}$ for any positive integers $n$ and $j$, which implies
$P(V_{j:n+j-1}\le x)$ is non-increasing in $j$ for any positive integer $n$ and $x\in (0,1)$. Since $\mbox{Beta}(j, l_r)$ and
$V_{j: l_r+j-1}$ have the same distribution, we have, for each $r$ and $x\in (0,1)$, $P(\mbox{Beta}(j, l_r)\le x)$ is non-increasing in
$j$. This concludes \eqref{stochastic-order}.


Now we prove \eqref{stability}.
Under condition $\lim_{n\to\infty}F_n(x)=F(x)$ for $x\in [0,1]$ together with (\ref{C_1_condition}) or (\ref{C_2_condition}), we first claim that
\begin{equation}\label{estimate-lj}
\frac{1}{\gamma_n}\sum^{m}_{r=1}\frac{l_r}{nx+l_j}=O(1)
\end{equation}
as $n\to\infty$ for any $x\in (0,1).$ First, from \eqref{eta-est} we have for any $\delta\in (0,1)$
\begin{equation}\label{logt}
-\frac{1+\delta}{2\delta}(1-t) \le \log t\le  -(1-t)~~~~~\mbox{  for } \delta\le t\le 1.
\end{equation}
Recall that $F_n(x)=\big(\prod^{m}_{r=1}\frac{n_rx}{nx+l_r}\big)^{1/\gamma_n}$ for $0< x<1$.
By assumption $\lim_{n\to\infty}F_n(x)=F(x)$ for $x\in [0,1]$,
\[
\frac{1}{\gamma_n}\sum^{m}_{r=1}\log\frac{n_rx}{nx+l_r}=O(1).
\]
Since $x\le \frac{n_rx}{nx+l_r}\le 1$ for any $1\le r\le m$, combining \eqref{logt} and the above equation we have that
for any $0<x<1$,
\[
\frac{1}{\gamma_n}\sum^{m}_{r=1}\frac{l_r(1-x)}{nx+l_r}=\frac{1}{\gamma_n}\sum^{m}_{r=1}
\Big(1-\frac{n_rx}{nx+l_r}\Big)=O(1),
\]
 which implies \eqref{estimate-lj} since $1-x>0$.

Note that
\begin{equation}\label{mean-var}
\mu_{j,r}=\E(s_{j,r})=\frac{j}{j+l_r},~~~\sigma^2_{j,r}=\Var(s_{j,r})=\frac{jl_r}{(j+l_r)^2(j+l_r+1)}.
\end{equation}
Recall (\ref{rep1}). From \eqref{beta-gamma} and \eqref{logY} we can rewrite the moment generating function of $\log Y_{n,j}^2$ as
\[
M_j(t)=\prod^m_{r=1}\Big(\frac{\Gamma(j+t)}{\Gamma(j)}\cdot\frac{\Gamma(j+l_r)}{\Gamma(j+t+l_r)}\Big).
\]
Hence we obtain
\[
\E(\log Y_{n,j}^2)=\frac{d}{dt}M_j(t)\Big|_{t=0}=\Big[M_j(t)\frac{d}{dt}(\log M_j(t))\Big]\Big|_{t=0}= \sum^m_{r=1}\Big(\frac{\Gamma'(j)}{\Gamma(j)}-\frac{\Gamma'(j+l_r)}{\Gamma(j+l_r)}\Big),
\]
that is,
\begin{equation}\label{logmean}
\E(\log Y_{n,j}^2)=\sum^m_{r=1}\big(\psi(j)-\psi(j+l_r)\big),
\end{equation}
where $\psi(x)=\Gamma'(x)/\Gamma(x)$ is the  Digamma function as mentioned in (\ref{ice_coffee}). Recall $\eta(x)=x-1-\log x$ for $x>0$, which is defined right before (\ref{eta-est}). Then
\begin{equation}\label{decomposition}
\log Y_{n,j}^2-\sum^{m}_{r=1}\log\mu_{j, r}=\sum^{m}_{r=1}\log \frac{s_{j,r}}{\mu_{j, r}}
=\sum^{m}_{r=1}\big(\frac{s_{j,r}}{\mu_{j, r}}-1\big)-\sum^{m}_{r=1}\eta\big(\frac{s_{j,r}}{\mu_{j,r}}\big).
\end{equation}
From now on, for each $x\in (0,1)$, we take $j=[nx]$. We will show
\begin{equation}\label{convergence1}
\frac{1}{\gamma_n}\sum^{m}_{r=1}(\frac{s_{j,r}}{\mu_{j,r}}-1)\overset{p}\to 0
\end{equation}
and
\begin{equation}\label{convergence2}
\frac{1}{\gamma_n}\sum^{m}_{r=1}\eta(\frac{s_{j,r}}{\mu_{j, r}}) \overset{p}\to 0.
\end{equation}

To prove \eqref{convergence1}, it suffices to show that the variance of the left-hand side in \eqref{convergence1} converges
to zero. In fact, we have from \eqref{mean-var} and \eqref{estimate-lj}
\begin{eqnarray*}
\Var\Big(\frac{1}{\gamma_n}\sum^{m}_{r=1}\big(\frac{s_{j,r}}{\mu_{j,r}}-1\big)\Big)
&=&\frac{1}{\gamma_n^2}\sum^{m}_{r=1}\frac{\sigma_{j, r}^2}{\mu_{j, r}^2}\\
&=&\frac{1}{\gamma_n^2}\sum^{m}_{r=1}\frac{l_r}{j(j+l_r+1)}\\
&\le&\frac{1}{\gamma_n^2}\sum^{m}_{r=1}\frac{l_r}{[nx]([nx]+l_r+1)}\\
&=&\frac{O(1)}{n\gamma_n} \sum^{m}_{r=1}\frac{l_r}{nx+l_r}\\
&\to&0
\end{eqnarray*}
as $n\to\infty$ by the assumption $\gamma_n\geq 1$.

Since $\frac{1}{\gamma_n}\sum^{m}_{r=1}\eta(\frac{s_{j,r}}{\mu_{j, r}})\ge 0$, to show \eqref{convergence2}, it suffices to verify that $\frac{1}{\gamma_n}\E\sum^{m}_{r=1}\eta(\frac{s_{j,r}}{\mu_{j, r}})\to 0$ as $n\to\infty$. To this end,
we have from \eqref{decomposition}, \eqref{mean-var} and \eqref{logmean} that
\begin{eqnarray*}
\E\Big[\sum^{m}_{r=1}\eta\big(\frac{s_{j,r}}{\mu_{j, r}}\big)\Big]&=&\Big(\sum^{m}_{r=1}\log\mu_{j, r}\Big)-\E\log( Y_{n,j}^2)\\
&=&\sum^{m}_{r=1}\Big[(\psi(j+l_r)-\log(j+l_r))-(\psi(j)-\log j)\Big]\\
&=&\sum^{m}_{r=1}\int^{j+l_r}_j(\psi'(t)-\frac{1}{t})dt.
\end{eqnarray*}
By formula 6.4.12 from Abramowitz and Stegun (1972)
\[
\psi'(t)=\frac{1}{t}+\frac{1}{2t^2}+O(\frac1{t^3})
\]
as $t\to+\infty.$ This and the fact $j=[nx]$ lead to
\begin{eqnarray*}
 \frac{1}{\gamma_n}\E\sum^{m}_{r=1}\eta\big(\frac{s_{j,r}}{\mu_{j, r}}\big)&=&\frac{1+o(1)}{\gamma_n}\sum^{m}_{r=1}\int^{j+l_r}_j\frac{1}{2t^2}dt\\
 &=& \frac{1+o(1)}{\gamma_n}\sum^{m}_{r=1}\frac12(\frac{1}{j}-\frac{1}{j+l_r})\\
 &=& \frac{1+o(1)}{\gamma_n}\sum^{m}_{r=1}\frac12\frac{l_r}{j(j+l_r)}\\
 &=& O(\frac{1}n)\frac{1}{\gamma_n}\sum^{m}_{r=1}\frac{l_r}{nx+l_r}\\
 &\to& 0
 \end{eqnarray*}
by \eqref{estimate-lj}.

From \eqref{decomposition} - \eqref{convergence2} we have for any $0<x<1$
\[
\frac{1}{\gamma_n}\Big(\log Y_{n,[nx]}^2-\sum^{m}_{r=1}\log\mu_{[nx],r}\Big)\overset{p}\to 0.
\]
Under (\ref{C_1_condition}) or (\ref{C_2_condition}), the limit $F(x)$  is continuous and positive in $(0,1)$. Therefore, the convergence $\lim_{n\to\infty}F_n(x)=F(x)$ is uniform for any interval $[\delta_1, \delta_2]\subset (0,1)$.  
It follows that  $\lim_{n\to\infty}F_n(\frac{[nx]}{n})=F(x)$ for any $x\in (0,1)$. Further,  notice that $b_n=g_n(1)$ and
$\sum^{m}_{r=1}\log\mu_{[nx],r}=\log g_n(\frac{[nx]}{n})$ by (\ref{mean-var}).
Then we have from (\ref{fnx}) that
\begin{eqnarray*}
\frac{1}{\gamma_n}\log \frac{Y_{n,[nx]}^2}{b_n}-\log F(x)
&=&\frac{1}{\gamma_n}\Big(\log Y_{n,[nx]}^2-\sum^{m}_{r=1}\log\mu_{[nx],r}\Big)+
\frac1{\gamma_n}\log\frac{g_n(\frac{[nx]}{n})}{g_n(1)}-\log F(x)\\
&=&\frac{1}{\gamma_n}\Big(\log Y_{n,[nx]}^2-\sum^{m}_{r=1}\log\mu_{[nx],r}\Big)
+\log F_n\Big(\frac{[nx]}{n}\Big)-\log F(x)
\end{eqnarray*}
converges in probability to zero for any $x\in (0,1).$ This completes the proof of \eqref{stability}. \hfill$\blacksquare$

\medskip

\noindent\textbf{Proof of Theorem \ref{pie}}. Easily, part (c) is a corollary of (a) and (b). So we only need to prove (a) and (b).

Notice, with the transform $Z=Re^{i\Theta}$, that the density of $Z$ is $\frac{1}{2\pi|z|}f^*(|z|)$ for $0<|z|< 1$ is equivalent to that the density of $(\Theta, R)$ is $\frac{1}{2\pi}f^*(r)$ for $\theta \in [0, 2\pi)$ and $0<r< 1$.

Let $\{s_{j,r},\,1\le j\le n,  1\leq r\leq m\}$ be independent random variables and $s_{j,r}$ have the Beta  density $\frac{1}{B(j,l_r)}y^{j-1}(1-y)^{l_r-1}I(0\leq y\leq 1)$ for each $j$ and $r$. By Lemma \ref{new_loving}, for ease of notation we assume, without loss of generality,  that
\begin{equation}\label{rep}
|Z_j|^2=\prod^m_{r=1}s_{j, r},  ~~~1\le j\le n.
\end{equation}

If condition (\ref{C_1_condition}) holds, since  $F_n(0)=0$, $F_n(1)=1$ and $\lim_{n\to\infty}F_n(x)=F(x)$ for $x\in (0,1)$, we assume without loss of generality that $\lim_{n\to\infty}F_n(x)=F(x)$ for $x\in [0,1]$, $F(0)=0$, $F(1)=1$, and $F(x)$ is continuous and strictly increasing over $[0,1]$. Hence,  $F^*(y):=F^{-1}(y)$ is a continuous and strictly increasing  distribution function on $[0, 1]$ with $F^*(0)=0$ and $F^*(1)=1$. Define $F^*(y)=1$ for $y>1.$ Further, $Unif\{|z|=1\}$ is identical to the product measure of $Unif[0, 2\pi)$ and $\delta_1$. The distribution function of $\delta_1$ is the indicator function $I_{[1, \infty)}(y)$. Let $Y_1, \cdots, Y_n$ be as in Lemma~\ref{independence} with $\varphi=w_m^{(l_1, \cdots, l_m)}$ defined in (\ref{truncation-density}). According to Theorem \ref{nonlinear} and Lemma~\ref{lem2}, to prove (a) and (b), it suffices to verify that
\bea\lbl{special_nude}
& & \lim_{n\to\infty}\frac{1}{n}\sum^n_{j=1}P\Big(\Big(\frac{Y_j^2}{b_n}\Big)^{1/\gamma_n}\le y\Big)\nonumber\\
& = &
\begin{cases}
F^*(y) & \text{for all $y>0$ if (\ref{C_1_condition}) holds;}\\
I_{[1, \infty)}(y) & \text{for all $y\in (0,1)\cup (1, \infty)$ if (\ref{C_2_condition}) holds.}
\end{cases}
\eea
By Lemma~\ref{independence},
\beaa
\sum^n_{j=1}P\Big(\Big(\frac{Y_j^2}{b_n}\Big)^{1/\gamma_n}\le y\Big)=\sum^n_{j=1}P\Big(\Big(\frac{|Z_j|^2}{b_n}\Big)^{1/\gamma_n}\le y\Big)
\eeaa
for any $y\in \mathbb{R}.$ Since $|Z_j|^2 \in [0, 1]$ is continuous for each $j$,
it suffices to show
\begin{equation}\label{cheers}
\lim_{n\to\infty}\frac{1}{n}\sum^n_{j=1}P\Big(\frac{1}{\gamma_n}\log \frac{|Z_j|^2}{b_n}\le \log y\Big)= \mbox{ the right hand side of (\ref{special_nude}). }
\end{equation}
From \eqref{stochastic-order} and \eqref{rep},
\bea\lbl{volleyball}
P\Big(\frac{1}{\gamma_n}\log \frac{|Z_j|^2}{b_n}\le \log y\Big)\ \mbox{is non-increasing in}\ j\in\{1, \cdots, n\}.
\eea
This property and equation \eqref{stability} play a central role in the following estimation.

\noindent{\it Proof of (\ref{cheers})  if} (\ref{C_2_condition}) holds.

Fix $y\in (0,1)$. 
For any $\delta\in
(0,1)$,  we have
\begin{eqnarray*}
&&\limsup_{n\to\infty}\frac{1}{n}\sum^n_{j=1}P\Big(\frac{1}{\gamma_n}\log \frac{|Z_j|^2}{b_n}\le \log y\Big)\\
&=&\limsup_{n\to\infty}\frac{1}{n}\Big[\sum^{[n\delta]}_{j=1}P\Big(\frac{1}{\gamma_n}\log \frac{|Z_j|^2}{b_n}\le \log y\Big)+\sum^{n}_{j=[n\delta]+1}P\Big(\frac{1}{\gamma_n}\log \frac{|Z_j|^2}{b_n}\le \log y\Big)\Big]\\
&\le&\limsup_{n\to\infty}\frac{[n\delta]}{n}+\limsup_{n\to\infty}P\Big(\frac{1}{\gamma_n}\log \frac{|Z_{[n\delta]}|^2}{b_n}\le \log y\Big)\\
&=&\delta,
\end{eqnarray*}
by the fact $\frac{1}{\gamma_n}\log\frac{Z_{[nx]}^2}{b_n}\overset{p}\to 0$ as
$n\to\infty$ for all $x\in (0,1)$ from (\ref{stability}) and (\ref{C_2_condition}). This implies
\bea\lbl{why_world}
\lim_{n\to\infty}\frac{1}{n}\sum^n_{j=1}P\Big(\frac{1}{\gamma_n}\log \frac{|Z_j|^2}{b_n}\le \log y\Big)=0
\eea
for $y\in (0,1).$
If $y>1$, then $\log y>0$. By (\ref{volleyball}),
\begin{eqnarray*}
1 &\geq &\liminf_{n\to\infty}\frac{1}{n}\sum^n_{j=1}P\Big(\frac{1}{\gamma_n}\log \frac{|Z_j|^2}{b_n}\le \log y\Big)\\
&\geq &\liminf_{n\to\infty}\frac{1}{n}\sum^{[n\tau]}_{j=1}P\Big(\frac{1}{\gamma_n}\log \frac{|Z_j|^2}{b_n}\le \log y\Big)\\
&\geq&\liminf_{n\to\infty}\frac{[n\tau]}{n}P\Big(\frac{1}{\gamma_n}\log \frac{|Z_{[n\tau]}|^2}{b_n}\le \log y\Big)=\tau
\end{eqnarray*}
as $n\to\infty$ for all $\tau \in (0,1)$. Letting $\tau \uparrow 1$, and combining with (\ref{why_world}), we get (\ref{cheers}).

\noindent{\it Proof of (\ref{cheers})  if} (\ref{C_1_condition}) holds. We will differentiate two cases: $y \in (0, 1)$ and $y\geq 1$.

\underline{Case 1: $y\in (0,1)$}. Set $F^*(y)=F^{-1}(y)$ for $y\in (0,1).$    Let $\delta\in (0,1)$ be a number such that  $0<y-\delta<y+\delta<1$.  Then
$0<F^*(y-\delta)<F^*(y+\delta)<1$, and $F(F^*(y-\delta))<y<F(F^*(y+\delta))$.  By letting $x=F^*(y+\delta)$, we have
\bea
&&\limsup_{n\to\infty}\frac{1}{n}\sum^n_{j=1}P\left(\frac{1}{\gamma_n}\log \frac{|Z_j|^2}{b_n}\le \log y\right)\nonumber\\
&=&\limsup_{n\to\infty}\Big[\frac{1}{n}\sum^{[nx]}_{j=1}P\Big(\frac{1}{\gamma_n}\log \frac{|Z_j|^2}{b_n}\le \log y\Big)+\frac{1}{n}\sum^{n}_{j=[nx]+1}P\Big(\frac{1}{\gamma_n}\log \frac{|Z_j|^2}{b_n}\le \log y\Big)\Big]\nonumber\\
&\le&\limsup_{n\to\infty}\frac{[nx]}{n}+\limsup_{n\to\infty}P\Big(\frac{1}{\gamma_n}\log \frac{|Z_{[nx]}|^2}{b_n}\le \log y\Big)\nonumber\\
&=&x+\limsup_{n\to\infty}P\Big(\frac{1}{\gamma_n}\log \frac{|Z_{[nx]}|^2}{b_n}-\log F(x)\le -\log\frac{F(x)}{y}\Big)\nonumber\\
&=&F^*(y+\delta) \lbl{fire_sign}
\eea
 by \eqref{stability} and the fact $\log(F(x)/y)>0$.   Similarly, setting $x=F^*(y-\delta)$, we get
 \bea
&&\liminf_{n\to\infty}\frac{1}{n}\sum^n_{j=1}P\Big(\frac{1}{\gamma_n}\log \frac{|Z_j|^2}{b_n}\le \log y\Big) \nonumber\\
&\ge &\liminf_{n\to\infty}\frac{1}{n}\sum^{[nx]}_{j=1}P\Big(\frac{1}{\gamma_n}\log \frac{|Z_j|^2}{b_n}\le \log y\Big) \nonumber\\
&\ge&\liminf_{n\to\infty}\frac{[nx]}{n}P\Big(\frac{1}{\gamma_n}\log \frac{|Z_{[nx]}|^2}{b_n}\le \log y\Big) \nonumber\\
&=&x\cdot\liminf_{n\to\infty}P\Big(\frac{1}{\gamma_n}\log \frac{|Z_{[nx]}|^2}{b_n}-\log F(x)\le \log\frac{y}{F(x)}\Big) \nonumber\\
&=&x\cdot \Big[1-\limsup_{n\to\infty}P\Big(\frac{1}{\gamma_n}\log \frac{|Z_{[nx]}|^2}{b_n}-\log F(x)>\log\frac{y}{F(x)}\Big)\Big] \nonumber\\
&=&F^*(y-\delta) \lbl{dinner_time}
\eea
by \eqref{stability} and the assertion $\log(y/F(x))>0$. Finally, by letting $\delta \downarrow 0$  in (\ref{fire_sign}) and (\ref{dinner_time}), we show (\ref{cheers}) holds under the condition  (\ref{C_1_condition}) and $y\in (0,1)$.

\underline{Case 2: $y\geq 1$}. Observe
\beaa
\frac{1}{n}\sum^n_{j=1}P\Big(\frac{1}{\gamma_n}\log \frac{|Z_j|^2}{b_n}\le \log y\Big) \geq \frac{1}{n}\sum^n_{j=1}P\Big(\frac{1}{\gamma_n}\log \frac{|Z_j|^2}{b_n}\le \log y_1\Big)
\eeaa
for all $y_1\in (0,1).$ By the proved conclusion, it is seen that
\beaa
\liminf_{n\to\infty}\frac{1}{n}\sum^n_{j=1}P\Big(\frac{1}{\gamma_n}\log \frac{|Z_j|^2}{b_n}\le \log y\Big) \geq F^*(y_1).
\eeaa
Then (\ref{cheers}) follows by taking $y_1\uparrow 1$ and by the fact $F^*(1)=1$. The proof is complete.    \hfill$\blacksquare$

\medskip

Now  we present the proofs of the corollaries.

\medskip

\noindent\textbf{Proof of Corollary \ref{cor1}}. 
Since $n_j=l_j+n$, take   $\gamma_n= 2$ for all $n\geq 1$ to have
\[
F_n(x)=\Big(\prod^{m}_{j=1}\frac{n_jx}{n_j-n(1-x)}\Big)^{1/2},  ~~~~0\le x\le 1.
\]
By assumption,  $\lim_{n\to\infty}\frac{n}{n_j}=\alpha_j\in [0,1]$ for $1\le j\le m$, we have
\begin{equation}\label{case1:F}
\lim_{n\to\infty}F_n(x)=\Big(\prod^{m}_{j=1}\frac{x}{1-\alpha_j(1-x)}\Big)^{1/2} :=F(x),\ \ \ 0<x\leq 1.
\end{equation}

\noindent (1). If $\alpha_1=\cdots=\alpha_m=1$, then $F(x)=1$ for $0< x\leq 1$.
That is,  condition (\ref{C_2_condition}) holds. The conclusion follows from (c) of Theorem \ref{pie}.

\medskip

\noindent (2) \& (3).
If $\alpha_j<1$ for some $1\leq j \leq m$, then $F(x)$ in (\ref{case1:F}) is continuous and strictly increasing in $(0,1]$. Also,  $\lim_{x\downarrow 0}F(x)=0$. The statement (3) then follows from (a) of Theorem \ref{pie}. In particular,
when $\alpha_1=\cdots=\alpha_m=\alpha\in [0,1)$, then $F(x)=\big(\frac{x}{1-\alpha(1-x)}\big)^{m/2}$ for $0<x\leq 1$. It is trivial to check that
\begin{equation*}
F^{-1}(x)=\frac{(1-\alpha)x^{2/m}}{1-\alpha x^{2/m}} ~~~~ \mbox{and} ~~~~ \frac{d}{dx}F^{-1}(x)=\frac{2(1-\alpha)}{m}\frac{x^{(2/m)-1}}{(1-\alpha x^{2/m})^2}.
\end{equation*}
We get (2). \hfill$\blacksquare$

\medskip

\noindent\textbf{Proof of Corollary \ref{cor2}}. 
Take $\gamma_n=m$. We will show that $\lim_{n\to\infty}F_n(x)=F(x)$ with
\bea\lbl{cailiflower}
F(x)=x\exp\Big(-\int^1_0\log(1-q(t)(1-x))dt\Big), ~~~0< x\le 1,
\eea
and $\lim_{x\downarrow 0}F(x)=0$. Obviously, $F(1)=1$. This says that condition (\ref{C_1_condition}) is satisfied.  We first prove (\ref{cailiflower}).   Since
\beaa
\log F_n(x) &=& \frac{1}{m}\sum^{m}_{j=1}\Big[\log x-\log \Big(1-\frac{n}{n_j}(1-x)\Big)\Big]\\
& = &
\log x- \frac{1}{m}\sum^{m}_{j=1}\log\Big(1-\frac{n}{n_j}(1-x)\Big),
\eeaa
we have that
\begin{eqnarray*}
|\log F_n(x)-\log F(x)|&=& \Big|\frac{1}{m}\sum^{m}_{j=1}\log\Big(1-\frac{n}{n_j}(1-x)\Big)-\int^1_0\log(1-q(t)(1-x))dt\Big|\\
&\le& \Big|\frac{1}{m}\sum^{m}_{j=1}\Big[\log(1-\frac{n}{n_j}(1-x))-\log\Big(1-q\Big(\frac{j}{m}\Big)(1-x)\Big)\Big]\Big|\\
&&~+\Big|\frac{1}{m}\sum^{m}_{j=1}\log\Big(1-q\Big(\frac{j}{m}\Big)(1-x)\Big)-\int^1_0\log(1-q(t)(1-x))dt\Big|\\
&:=&I_1(x)+I_2(x)
\end{eqnarray*}
for $0<x\leq 1$. Since for each fixed $x\in (0,1]$,  $\frac{d}{dt}\log(1-t(1-x))=\frac{1-x}{1-t(1-x)}$ for $0\le t\le 1$,  we have
$1-x\le \frac{d}{dt}\log(1-t(1-x))\le \frac{1-x}{x}$. By using the mean-value theorem and (\ref{case22}), we have
\[
I_1(x)\le \frac{1-x}{x}\frac{1}{m}\sum^{m}_{j=1}\Big|\frac{n}{n_j}-q\Big(\frac{j}{m}\Big)\Big|\to 0
\]
as $n\to\infty$. For fixed $x\in (0,1]$,  since $\log(1-t(1-x))$ is a bounded and continuous function in $t\in [0,1]$ and thus
$\log(1-q(t)(1-x))$ is also a bounded and continuous function in $t\in [0,1]$, we have $I_2(x)\to 0$ as $n\to\infty$ by the definition of the Riemann integral. Therefore, we have proved that $\lim_{n\to\infty}F_n(x)=F(x)$ for $x\in (0,1]$.

Choose any $t_0\in (0,1)$ such that $q(t_0)<1.$ From continuity, there exists $\epsilon>0$ satisfying $0< t_0-\epsilon< t_0+\epsilon<1$ and $\sup_{|t-t_0|\leq \epsilon}q(t)\leq 1-\epsilon$. Easily, $\frac{x}{1-q(t)(1-x)}< 1$ for all $x\in (0,1)$. It follows that
\begin{eqnarray*}
\log F(x)&=&\log x-\int^1_0\log(1-q(t)(1-x))dt\\
&=&\int^1_0\log\frac{x}{1-q(t)(1-x)}dt\\
&\le& \int^{t_0+\epsilon}_{t_0-\epsilon}\log\frac{x}{1-q(t)(1-x)}dt\\
&\le& \log\frac{2\epsilon x}{1-(1-\epsilon)(1-x)}\\
&\to &-\infty
\end{eqnarray*}
as $x\downarrow 0$, yielding that $\lim_{x\downarrow 0}F(x)=0$. Thus, (\ref{C_1_condition}) is verified.

Now,
\beaa
f(x)=F'(x) & = & \frac{1}{x}F(x) - F(x)\int_0^1\frac{q(t)}{1-q(t)(1-x)}\,dt\\
& = & F(x)\int_0^1\frac{1-q(t)}{1-q(t)(1-x)}\,dt.
\eeaa
Further,
\[
f^*(x)=\frac{d}{dx}F^{-1}(x)=\frac{1}{f(F^{-1}(x))},  ~~~0<x<1.
\]
By (a) of Theorem \ref{pie}, $\mu_n \rightsquigarrow\mu$, where $\mu$ has density $\frac{1}{2\pi}f^*(r)$ for $(\theta, r)\in [0, 2\pi)\times [0,1)$.\hfill$\blacksquare$

\medskip

\noindent\textbf{Proof of Corollary \ref{cor3}}. 
Fix $0<x\leq 1$. Then,
\bea\lbl{weapon}
\log F_n(x)=-\frac{1}{\gamma_n}\sum_{j=1}^{m}\log \Big[\frac{1}{x}\Big(1+\frac{n}{n_j}(x-1)\Big)\Big].
\eea
Write
\beaa
\frac{1}{x}\Big(1+\frac{n}{n_j}(x-1)\Big) =1 + \frac{1-x}{x}\Big(1-\frac{n}{n_j}\Big).
\eeaa
By the assumption $\lim_{n\to\infty}\max_{1\le j\le m}|\frac{n}{n_j}-1|=0$, we know the logarithm on the right side of (\ref{weapon}) is equal to $\frac{1-x}{x}\big(1-\frac{n}{n_j}\big)(1+o(1))$ uniformly for all $1\leq j \leq m$. Use the identity that $n_j=l_j+n$ to see
\bea\lbl{sea_mirror}
\log F_n(x) &=& \frac{x-1}{\gamma_n\,x}(1+o(1))\sum_{j=1}^{m}\frac{l_j}{n_j} \nonumber\\
& = & (1+o(1))\frac{x-1}{x}\cdot\frac{1}{\gamma_n}\cdot \frac{1}{n}\sum_{j=1}^{m}l_j,\ \ \ 0<x\leq 1,
\eea
as $n\to\infty$ by the assumption  $\lim_{n\to\infty}\max_{1\le j\le m}|\frac{n}{n_j}-1|=0$ again. We now show (a), (b) and (c) according to the assumption that $\lim_{n\to\infty}\frac1n\sum^{m}_{j=1}l_j=\beta\in [0,\infty]$.

\noindent(a). Assume $\beta=0$. Choose $\gamma_n=2$. Then, $F(x)=1$ for $0<x\leq 1.$ The conclusion follows from the second part of (c) from Theorem~\ref{pie}.

\noindent(b). Assume $\beta \in (0, \infty)$. Choose $\gamma_n=2$. Then $F(x)=\exp(\frac{\beta}{2}\frac{x-1}{x})$ for $0<x\leq 1$ and $\lim_{x\downarrow 0}F(x)=0$. Trivially, $F^{-1}(x)=\big(1-\frac{2}{\beta}\log x\big)^{-1}$ for $0<x \leq 1$ and
\beaa
\frac{d}{dx}F^{-1}(x)=\frac{2\beta}{x(\beta-2\log x)^2}, \ \ \ 0<x \leq 1.
\eeaa
So the density of $\mu$ according to (a) of Theorem~\ref{pie} is
\beaa
f(\theta,r):=\frac{\beta}{\pi r(\beta-2\log r)^2}
\eeaa
for $(\theta, r)\in [0, 2\pi)\times [0, 1).$ By the polar transformation $z=r e^{i\theta}$, it is easy to see that
the density of $\mu$ is given by
 \beaa
 h(z):=\frac{\beta}{\pi |z|^2(\beta-2\log |z|)^2},\ \ |z|<1.
 \eeaa

\noindent(c). Assume $\beta=\infty$. Take $\gamma_n=\frac1n\sum^{m}_{j=1}l_j.$ Then, by (\ref{sea_mirror}), $F(x)=\exp\{\frac{x-1}{x}\}$ for $0< x\leq 1$ and $\lim_{x\downarrow 0}F(x)=0$. This is the same as the case that $\beta=2$ in (b). So the density of $\mu$ is $\frac{1}{2\pi |z|^2(1-\log |z|)^2}$ for $0<|z|<1.$    \hfill$\blacksquare$

\medskip

\noindent\textbf{Proof of Corollary \ref{cor4}}. By assumption, $m\to\infty$  and $\max_{1\le j\le m}\frac{n}{n_j}=0$ as $n\to\infty$. By using the equality $n_j=l_j+n$ and taking $\gamma_n=m$, we see that
\[
\lim_{n\to\infty}F_n(x)=x:=F(x)
\]
for any $x\in [0,1]$. Then $F^{-1}(x)=x$ for $x\in [0, 1]$ and $f^*(x)=1$ for $0\leq x \leq 1$.
Choosing $h_n(r)=(r^2/b_n)^{1/\gamma_n}$ and  $b_n=\prod^{m}_{j=1}\frac{n}{n_j}$, we know that $\mu_n  \rightsquigarrow (\Theta, R)$ and $Z=R e^{i\Theta}$ has the  density $\frac{1}{2\pi |z|}$ for $0<|z|\leq 1$. Equivalently, $(\Theta, R)$ has law $Unif([0, 2\pi)\times [0,1]).$ As mentioned before, $Z=Re^{i\Theta}$ has the uniform distribution on $\{z\in \mathbb{C};\, |z|\leq 1\}$ if and only if $(\Theta, |Z|^2)$ follows the uniform distribution on $[0, 2\pi)\times [0,1).$  Therefore,  the desired conclusion follows. \hfill$\blacksquare$\\

 \noindent\textbf{Acknowledgements}. We thank an anonymous referee for his/her very careful reading. The referee's report helps us make the presentation much more clearly.

\newpage

\baselineskip 12pt
\def\ref{\par\noindent\hangindent 25pt}

\end{document}